\documentclass[12pt]{amsart}
\usepackage{amsmath,amsthm,amssymb,xypic}

\theoremstyle{plain}  

\newtheorem{thm}{Theorem}[section]

\newtheorem{lem}[thm]{Lemma}   
\newtheorem{cor}[thm]{Corollary}
\newtheorem{con}[thm]{Conjecture}

\newtheorem{rem}[thm]{Remark}

\theoremstyle{definition}

\newtheorem{dfn}[thm]{Definition}

\theoremstyle{remark}

\newtheorem{exe}[thm]{Exercise}

\DeclareMathOperator{\Cl}{Cl}

\DeclareMathOperator{\NEbar}{\overline{NE}}

\DeclareMathOperator{\bs}{Bs}

\DeclareMathOperator{\hcf}{hcf}

\DeclareMathOperator{\mult}{mult}

\DeclareMathOperator{\Exc}{Exc}

\DeclareMathOperator{\Proj}{Proj}

\DeclareMathOperator{\Sing}{Sing}
\DeclareMathOperator{\NonSing}{NonSing}
\DeclareMathOperator{\WCl}{WCl}
\DeclareMathOperator{\Pic}{Pic}

\newcommand{\QED}{\ifhmode\unskip\nobreak\fi\quad {\rm Q.E.D.}} 

\newcommand\Ga{\Gamma}
\newcommand\flip{\dasharrow}

\newcommand\f{\varphi}
\newcommand\la{\lambda}
\newcommand{\wave}{\widetilde}
\newcommand{\fie}{\varphi}
  
\newcommand{\C}{\mathbb{C}}
\newcommand{\Z}{\mathbb{Z}}

\newcommand{\Q}{\mathbb{Q}}
\newcommand{\T}{\mathcal{T}}
\renewcommand{\P}{\mathbb{P}}
\renewcommand{\H}{\mathcal{H}} 
\renewcommand{\O}{\mathcal{O}}  
\renewcommand{\L}{\mathcal{L}}

\begin{document}

\title[Birational geometry of terminal quartic 3-folds. I]{Birational geometry of \\
terminal quartic 3-folds. I}

\author{
Alessio Corti
\and
Massimiliano Mella}
\address{A. Corti\\ 
D.P.M.M.S., Centre for Mathematical Sciences\\
University of Cambridge\\
Wilberforce Road\\ 
Cambridge CB3 0WB, U.K.}
\email{a.corti@dpmms.cam.ac.uk}
\address{M. Mella\\
Dipartimento di Matematica\\ 
Universit\`a di Ferrara\\
Via Machiavelli 35\\
44100 Ferrara Italia}
\email{mll@unife.it}
\date{January 2001}

\maketitle

\section{Introduction} \label{sec:intro}

\subsection{Abstract} 

In this paper, we study the birational geometry of 
certain examples of mildly singular quartic 3-folds. A quartic 3-fold
is a special case of a Fano 
variety, that is, a variety $X$ with ample anticanonical sheaf
$\O_X(-K_X)$. Nonsingular Fano 3-folds have been studied quite extensively.
From the point of view of birational geometry they basically fall within two 
classes: either $X$ is ``close to being rational'', and then it has very many 
biregularly distinct birational models as a Fano 3-fold, or, at the other 
extreme, $X$ has a unique model and it is often even true that every 
birational selfmap of $X$ is biregular. In this paper we construct examples
of singular quartic 3-folds with {\em exactly two} birational models as
Fano 3-folds; the other model is a complete intersection $Y_{3,4}
\subset \P(1,1,1,1,2,2)$ of a quartic and a cubic in weighted projective
space $\P(1^4, 2^2)$. These are the first examples to show this type
of behavior. 
After a brief introduction to singularities and Fano 3-folds, we give a first
precise statement of our main theorem. 
The rest of the introduction is not logically necessary to understand
the results or their statements. 
We describe the more general context of Fano 3-folds and Mori fiber 
spaces, and we outline a program of research on 3-folds which brings together
birational geometry, classification theory and commutative algebra.

\subsection{Singular quartic 3-folds}

We study quartic 3-folds having a unique singular point $P \in X$ 
analytically equivalent to $xy+z^3+t^3=0$. Choosing coordinates 
$x_0, ... x_4$ in $\P^4$ such that $P$ is the point $(1,0,...,0)$, 
we may write the equation of $X$ as
 \[
 F=x_0^2 x_1x_2+x_0a_3+b_4
 \]
where $a_3$, $b_4$ are homogeneous forms of the indicated degree in the 
variables $x_1, ...,x_4$. We also assume that $a_3$, $b_4$ are
general, where we take ``general'' in the sense of ``outside a Zariski
closed set'' in moduli. Two properties of $X$ are crucial to us.

The first is that the singularity $xy+z^3+t^3=0$ is {\em terminal}.
We need not enter into the precise details of the definition; the point
of terminal singularities is that certain manipulations with the canonical
class and discrepancy, familiar from the nonsingular case, still hold.
The most important examples are isolated hypersurface singularities of 
the form $xy+f(z,t)=0$, and quotients of $\C^3$ by the diagonal action of
$\Z/r\Z$ with weights $(1, a, -a)$ (when $a$ is coprime with $r$). The 
reader can look into \cite{YPG} for an accessible introduction to terminal 
singularities. 

The second property is that $X$ is $\Q$-factorial. This is an important and 
subtle (and much misunderstood) concept of Mori theory. It means by definition
that every Weil divisor of $X$ is $\Q$-Cartier; this is a property of the
Zariski topology of $X$ and not of the analytic type of its singularities.
If $X$ has hypersurface terminal singularities, $X$ is $\Q$-factorial 
if and only if it is factorial, that is all its (Zariski) local rings
are UFDs. In the case of a Fano 3-fold $X$, it is easy to see that 
$\Q$-factorial is equivalent to $\dim H^2(X)=\dim H_4(X)$, i.e., in
the case of a quartic
3-fold with terminal singularities, the 4th integral homology group 
$H_4(X, \Z)$ is generated by the class of a hyperplane section. This property
is often tricky to verify in practice; see below for additional comments on 
this. For a very gentle introduction to Mori theory, we recommend the 
Foreword to \cite{CR}.

\subsection{Fano 3-folds}

A quartic 3-fold is a particular case of a Fano 3-fold. There are 16 
deformation families of nonsingular Fano 3-folds with $\dim H^2=1$; with 
only one exception they were known classically.
Most of these varieties can be shown to be rational, and sometimes 
unirational, by means of classical constructions. 

Mori theory requires that we look at Fano 3-folds with terminal singularities
with $\dim H^2 =1$ and $\Q$-factorial; this is indeed one of the
possible end products of the minimal model program. In this paper, we
use the terminology
``Fano 3-fold'' in this more general sense implied by Mori theory.

There is at present no complete classification of Fano 3-folds, but 
several hundred families are known, see \cite{IF} and Alt{\i}nok \cite{Al} 
for some lists. For instance, there are 95 families of Fano 3-fold weighted 
hypersurfaces; it is conjectured in \cite{CPR}, and proved in many cases, 
that every quasismooth member of one of these families is the unique model
as a Fano 3-fold in its birational class.
There are also $85$ families of Fano 3-fold
codimension $2$ weighted complete intersections of which $Y_{3,4}$ is an 
example \cite{IF}, $69$ codimension $3$ Pfaffians and more than 
$100$ codimension $4$ families \cite{Al}.

Not only there are many more singular Fano 3-folds than nonsingular ones; 
they have a much richer variety of behavior from the point of view of 
birational geometry. In this paper we begin to study some of the simplest
new phenomena.

\subsection{Main result, crude statement}

This is our main result.

 \begin{thm} \label{thm:main} Let $X=X_4 \subset \P^4$ be a quartic 3-fold, 
with a singularity $P\in X$ analytically equivalent to $xy+z^3+t^3=0$, 
but otherwise general (in particular, nonsingular outside $P$). Then 

 \begin{enumerate}
 \item Let $Y$ be a Fano 3-fold (according to our conventions, 
this includes that $Y$ has terminal singularities, is $\Q$-factorial and 
has $\dim H^2X=1$). If $Y$ is birational to $X$, then {\em either}:
  \begin{enumerate}
  \renewcommand{\labelenumi}{(\alph{enumi})}
  \item $Y\cong X$ is biregular to $X$, {\em or}
  \item $Y=Y_{3,4}\subset \P(1,1,1,1,2,2)$ is a quasismooth complete 
  intersection of a quartic and a cubic in weighted projective space 
  $\P(1^4, 2^2)$.
  \end{enumerate}
 \item Let $Y$ be a 3-fold admitting either a morphism 
$Y\to S$ to a curve $S$ with typical fibre a rational surface, or 
$Y\to S$ to a surface $S$ with typical fibre $\P^1$. Then $Y$ is not
birational to $X$.
 \end{enumerate}
 \end{thm}

\begin{rem}
We could make the generality requirements on $X$ explicit if we
wanted; the precise conditions are $2.2(a-b)$ and the condition $(a)$
at the end of \S~6 (and~7). The result may well be true without the generality assumption, but it
would take more work to establish by the methods of this paper.

Note that $X$ is indeed factorial: even though the \textit{analytic} 
singularity $xy+z^3+t^3=0$ has a nontrivial class group (it is equivalent
to $xy+tz(t+z)=0$, so that $x=t=0$, for example, is not locally the divisor 
of a function), nevertheless when it appears on $X$
it is \textit{algebraically} factorial. 
We briefly explain why
   (this proof was suggested by J. Koll\'ar). The
   hyperplane section $X_0$ at infinity is generic if $X$ is generic,
   therefore by Lefschetz it has Picard rank one. The usual argument
   with formal schemes shows that $\Pic X \setminus P$ injects in
   $\Pic X_0$. This shows that a generic $X$ is $\Q$-factorial. But we
   know that $\Q$-factorial is a global topological property, and all
   $X$ having a unique singular point of the specified type are
   diffeomorphic. 
\end{rem}

In Section~\ref{sec:sbm} we exhibit an explicit birational 
map of $X_4$ to a quasismooth $Y_{3,4}\subset \P(1^4,2^2)$, the complete 
intersection of a cubic and a quartic in the weighted 5-dimensional 
projective space $\P(1,1,1,1,2,2)$. A striking feature of
this construction is that $Y_{3,4}$ is general; in fact any quasismooth 
$Y_{3,4}$ is birational to a quartic $X_4$ of the special kind considered.

In the rest of the introduction, we outline a program of research in
birational geometry generalizing the results of this paper in various
directions. 

\subsection{Birational geometry, classification theory and commutative algebra}

It is not difficult, and fun, for someone experienced in the use of 
the known lists of Fano 3-folds, and aware of certain constructions of 
graded rings, to generate many examples of birational maps between 
Fano 3-folds. In particular many Fano 3-fold codimension 2 weighted complete 
intersections are birational to special Fano 3-fold hypersurfaces, and we 
conjecture that a statement analogous to our main theorem holds for 
a lot of them. To name but a few, a general $Y_{6,7}\subset \P(1,1,2,3,3,4)$ is
birational to a special $X_7 \subset \P(1,1,1,2,3)$ with a singular
point $y^2+z^2+x_1^6+x_2^6$, a general $Y_{14, 15} \subset \P(1,2,5,6,7,9)$ 
is birational to a special $X_{15} \subset \P(1,1,2,5,7)$ with a singular
point $u^2+z^2y+y^7+x^{14}$, etc. It is remarkable that a significant part of
the list of Fano weighted complete intersections can be generated in this way,
starting from singular hypersurfaces.

This is a fairly general phenomenon. When trying to classify Fano
3-folds, the problem is often to construct a variety $Y$  with a given Hilbert 
function. Usually $Y$ has high codimension; in the absence of a structure
theory of Gorenstein rings, one method to construct $Y$ starts by studying
a suitable projection $Y \dasharrow X$ to a Fano $X$ in smaller codimension 
(the work of Fano, and then Iskovskikh, is an example of this). 
The classification of Fano 3-folds involves the study of the geometry of 
special members of some families (like our special singular
quartics), as well as general members of more complicated families 
(like our codimension $2$ complete intersections); the two points of view
match like the pieces of a gigantic jigsaw puzzle.

Miles Reid calls the map $X \dasharrow Y$ an ``unprojection''. Given $X$
satisfying certain properties, the problem is to construct $Y$. There are
at present a handful known constructions of this kind, leading to 
formats for codimension 4 Gorenstein rings which are often
a good substitute for the still missing general structure theory. 

The ideas here are due to Miles Reid, see for example \cite{k3}; for these 
and other issues not touched upon in this introduction, we also refer to 
\cite{r2}.

There is here an interplay of different problem areas: methods
of unprojection were first discovered in birational geometry   
\cite{CPR} 7.3, then applied to the construction of Fano 3-folds
\cite{Al}, and formats of Gorenstein rings \cite{Pa}. Commutative algebra
in turn clarifies birational geometry.

\subsection{The Sarkisov category}

 \begin{dfn} 

 \begin{enumerate}
 \item The {\em Sarkisov category} is the category
whose objects are Mori fibre spaces and whose morphisms birational maps
(regardless of the fibre structure).

 \item Let $X\to S$ and $X'\to S'$ be Mori fibre spaces.
A morphism in the Sarkisov category, that is, a birational map
$f\colon X\dasharrow X'$, is {\em square\/} if it
fits into a commutative square
 \[
 \mbox{\definemorphism{birto} \dashed \tip \notip
\diagram
X \rbirto^f \dto & X' \dto \\
S \rbirto^g & S' \enddiagram}
 \]
where $g$ is a birational map (which thus identifies the function field
$L$ of $S$ with that of $S'$) and, in addition, the induced birational
map of generic fibers $f_L\colon X_L\dasharrow X'_L$ is biregular. In this
case, we say that $X\to S$ and $X'\to S'$ are {\em square
birational}, or {\em square equivalent}.

 \item A {\em Sarkisov isomorphism} is a birational map
$f\colon X\dasharrow X'$ which is biregular and square.

 \item If $X$ is an algebraic variety, we define the {\em pliability} of
$X$ to be the set
 \[
\mathcal{P}(X)=\bigl\{ \text{Mfs} \; Y\to T\mid Y\;\text{is birational to} 
\;X \bigr\}/\text{square equivalence}.
 \]
We say that $X$ is {\em birationally rigid} if $\mathcal{P}(X)$ consists of 
one element.
 \end{enumerate}
 \end{dfn}

To say that $X$ is rigid means that it has an essentially unique model as 
a Mori fibre space; this implies in particular that $X$ is nonrational, but 
it is much more precise than that. For example, it is known that a nonsingular 
quartic 3-fold is birationally rigid \cite{IM}, \cite{Pu2}, \cite{Co2}; 
a quartic 3-fold which is nonsingular except for a single ordinary node is 
also birationally rigid \cite{Pu}, \cite{Co2}.

\subsection{The main result restated}

We restate our main theorem in the language of the Sarkisov category:

 \begin{thm} Let $X=X_4 \subset \P^4$ be a quartic 3-fold, 
with a singularity $P\in X$ analytically equivalent to $xy+z^3+t^3=0$, 
but otherwise general (in particular $X$ is nonsingular outside $P$). Then 
$\mathcal{P}(X)$ consists of two elements.
 \end{thm}

\subsection{Strict Mori fibre spaces}

Our focus in this paper is Fano 3-folds. It is natural, and eventually
necessary, to study similar problems in the more general context 
of strict Mori fibre spaces. For example Grinenko \cite{Gr} looks
at a ``double singular quadric'', i.e.,  $Z$ is a special $Z_{2,4}
\subset \P(1^5,2)$, where the degree $2$ equation is
the cone $x_1x_2+x_3x_4=0$ over a 2-dimensional quadric. Here $Z$
has two fibre structures of Del Pezzo surfaces of degree $2$, corresponding
to the rulings of the (2-dimensional) quadric, and it is shown that 
$\mathcal{P} (Z)$ consists of these two Mfs.

Similar examples are a general $Y_{3,3} \subset \P(1^5, 2)$ birational
to a cubic Del Pezzo fibration, but also a general
$Y_{4,4} \subset \P(1^3, 2^3)$ birational to a Del Pezzo fibration of degree
$2$, a general $Y_{6,6} \subset \P(1,1,2,2,3,4)$ birational to a Del Pezzo
fibration of degree $1$, etc. We may expect that these and many 
similar examples will be studied extensively in the near future.

\subsection{Pliability and rationality}

Traditionally, we like to think of Fano 3-folds as 
being ``close to rational''. We are now confronted with a view of 
3-fold birational geometry of great richness, on a scale much
larger than accessible with the methods of calculation and 
theoretical framework prior to Mori theory.

The notion of pliability is more flexible; a case division in terms of the 
various possibilities for $\mathcal{P}(X)$ allows to individuate a wider
spectrum of behavior ranging from birationally rigid to rational.

\subsection{Our starting point}

Our starting point and eventual goal is a uniform study of the
pliability of (singular) quartic 3-folds. We hope that we will soon be
able to settle the following

 \begin{con} Let $X=X_4\subset \P^4$ be a quartic 3-fold satisfying
the following conditions

 \begin{enumerate}
 \item $X$ has isolated singular points, that are locally analytically of 
the form $x^2+y^2+z^2+t^n=0$, for some positive integer $n$ (dependent
on the point in question),

 \item $X$ is factorial (this is equivalent to $\mathbb{Q}$-factorial,
   $X$ is a hypersurface).
 \end{enumerate}
Then $X$ is birationally rigid.
 \end{con}

\begin{rem}
According to Arnold, the next more complicated singularity is 
$xy+z^3+t^3=0$, and this is the case which we study here.
\end{rem}

Note the most striking feature of this conjecture: we do not 
restrict the number of singular points, although we insist that $X$ must be
factorial. We have seen that this is a subtle condition equivalent
to $\dim H_4=1$ i.e. the integral homology group $H_4(X, \Z)$ is generated by
the class of a hyperplane section. 

For example, if $X$ has an ordinary node as its unique singular point, then 
$X$ is automatically factorial. On the other hand, if $X$ has only ordinary 
nodes as singularities, $X$ is $\Q$-factorial if and only if the nodes impose
independent linear conditions on cubics. Indeed if $\tilde X$ is the blowup 
of the nodes, $H^1(\tilde X, \Omega_{\tilde X}^2)$ arises from residuation 
of 3-forms 
 \[
 \frac{P\Omega}{F^2}
 \]
on $\P^4$ with a pole of order two along $X= \{F=0\}$. Here $\Omega =
\sum x_i \, dx_0  \cdots \widehat{dx_i}\cdots dx_4$ and $P$ is a cubic
containing all the nodes of $X$, see e.g. \cite{Cl}. 

Consider a quartic 3-fold $Z$ containing the plane 
$x_0=x_1=0$. The equation of $Z$ can be written in the form 
$x_0a_3+x_1b_3=0$ and, in general, $Z$ has $9$ ordinary nodes 
$x_0=x_1=a_3=b_3=0$. The linear system $|a_3, b_3|$ defines a map
to $\P^1$; blowing up the base locus gives a Mori fibre space $Z \to
\P^1$ with fibers cubic Del Pezzo surfaces. However, our conjecture
does not apply to $Z$. Indeed $Z$ is not $\Q$-factorial:
the plane $\{x_0=x_1=0\}\subset Z$ is not a Cartier divisor. Thus 
$Z$ is not a Mori fibre space; it doesn't even make sense to say that it 
is rigid. (Note in passing: introducing the ratio $y=a_3/x_1=b_3/x_0$, 
gives a birational map $Z\dasharrow Y_{3,3}\subset \P(1^5, 2)$ to a
Fano 3-fold $Y_{3,3}$, the complete intersection of two cubics
in $\P(1^5, 2)$, a Mori fibre space birational to $Z$.
In the language of the Sarkisov program, $Z$ is the \textit{midpoint}
of a \textit{link} $X \dasharrow Y_{3,3}$.)

However, a quartic 3-fold with 9 nodes is factorial in general, and our 
conjecture predicts that then it is birationally rigid.

The factoriality of projective hypersurfaces is the subject of a
lovely paper by C. Ciliberto and V. Di Gennaro \cite{CDG}.

\subsection{Pliability and deformations}

We hope that the notion of pliability will be helpful in other ways too.
It is not known how rationality and, especially, unirationality behave in 
families. It is suspected that rationality is not stable under deformations.
For example, a general nonsingular quartic 4-fold $X$ in $\P^5$ is probably
nonrational, whereas $X$ is rational if it contains two skew planes.
The situation is even worse with unirationality, because there is no known 
way (other than trivial reasons) to show that a variety is not 
unirational; every judgment about the behavior
of this notion is at present pure guesswork. We can hope that 
$\mathcal{P}(X)$ is often a finite union of algebraic varieties, and that it 
is reasonably well behaved under deformations, for example, that 
it is an upper semi-continuous function of $X$.
For example, consider again a quartic $X$ containing a plane; 
then $X$ has $9$ nodes but it is not factorial. It is impossible
to deform $X$ to a quartic $X'$ with $9$ nodes in general position
(by what has been said, $\dim H_4X=2$ and $\dim H_4X'=1$, hence $X$ is not 
diffeomorphic to $X'$). This seems to suggest that birationally
rigid is quite robust under deformations.

The crucial point here is not so much whether this is literally true or not;
indeed it would be foolish to try to legislate over a large body of 
still largely unexplored examples. The point is that we now have the 
technology to test these ideas on substantial examples.

\subsection{Mori theory grows down}

Mori theory has enjoyed an initial phase of tremendous abstract development,
which continues today with higher dimensional flips and abundance.
Our main interest instead lies in the program of explicit birational 
geometry of 3-folds as described in the Foreword to \cite{CR}. The aim
is to treat 3-folds as explicitly as possible. An important
step in this development is to work out explicit description for the steps
of the minimal model program, divisorial contractions and flips, and the links
of the Sarkisov program.

\subsection{Extremal divisorial contractions}

An extremal divisorial contraction is a birational extremal contraction
$f\colon E\subset Z \to P \in X$, in the Mori category, which contracts a 
divisor $E$. This is the 3-fold analogue of the contraction of a $-1$ 
curve on a
surface. Despite some remarkable recent progress in special cases \cite{Kw}, 
we still don't know an explicit classification of 3-fold extremal divisorial
contractions. We prove in Section~\ref{sec:dco} that, if $P\in X$ is
the singularity $xy+z^3+t^3=0$, and $f$ contracts $E$ to the origin, 
then $f$ is the weighted blowup with weights $(2,1,1,1)$ or $(1,2,1,1)$. 
This classification is a key point in the proof of the main theorem.
In the proof we combine, and refine, two methods which have been
instrumental in the solution of similar problems. One is the
connectedness theorem of Shokurov, which is the key ingredient in
the classification of divisorial contractions in the case $P\in X$ 
is the ordinary node $xy+zt=0$ \cite{Co2}. The other
is a fairly simple minded multiplicity calculation which is sufficient
to establish the case $P\in X$ a terminal quotient singularity 
$1/r(a,-a,1)$ \cite{Ka}. Our proof is not easy and presupposes a good 
understanding of these other (simpler) cases. We advise the reader who is not 
an expert to study the simpler cases first, or else just skim 
through the proof.

\subsection{Methods}

Our proof uses the Sarkisov program \cite{Co1} and builds on and refines
the methods of \cite{CPR}, \cite{Co2}. 
The refinements are quite subtle and we expect that the nonspecialist
reader will find it difficult to follow the details of the 
proof. Because of this, we have made an effort to keep most of the discussion 
accessible to a general audience. 
We hope to have successfully swept the most technical parts of 
the proof under the final two Sections~\ref{sec:x4} and~\ref{sec:x34}.
These demand a great deal of motivation, and experience, on the part
of the reader.

In short, the new elements are 
 \begin{itemize}
 \item the partial classification of divisorial contractions, which improves
 on previous work, as already explained;
 \item a use of test surfaces and especially inequalities arising 
 from the theory of log surfaces, in combination with Shokurov's
 ``inversion of adjunction'', more efficient than for example in the proof of 
 rigidity of a general $Y_{2,3}\subset \P^5$ given in \cite{Co2};
 \item the general organization of the exclusion of curves as centers in 
 Section~\ref{sec:x34}, which circumvents the need of having to treat a large 
 number of particular cases.  
 \end{itemize}

\subsection{Organization and contents}

The paper is organized as follows. In Section~\ref{sec:sbm} we construct the 
birational map $X_4 \dasharrow Y_{3,4}$. We explain the construction at 
length within a narrative context and in much greater generality than needed
for the treatment of our example. Our aim is to equip the reader
with the knowledge to do many calculations of this type.

Section~\ref{sec:dco} is a short survey of what is known about the 
classification of 3-fold divisorial contractions. The material here is 
interesting in its own right; we tried to be as self-contained as possible, 
except for the proof of our own Theorem~\ref{thm:dco}.

In Section~\ref{sec:plan} we outline the proof of the main theorem. This is
copied from \cite{CPR} 3.1, 3.2, 3.3, with very few words changed. 

In Section~\ref{sec:excl} we survey the technique for excluding maximal 
centers, and indicate our own improvements for later use; in the final 
Sections~\ref{sec:x4} and \ref{sec:x34} we get down to the gruesome details of
the classification of  maximal singularities on $X_4$ and $Y_{3,4}$, namely 
we show that, other than those untwisted in Section \ref{sec:sbm}, there 
are no other maximal singularities. The results here combine to give a 
proof of the main theorem.

\subsection{Acknowledgements} We like to thank Miles Reid for following 
this project with interest and a huge amount of advise and teaching. 
In particular he helped us to understand the algebra of our map 
$X_4 \dasharrow Y_{3,4}$. We also like to thank J\'anos Koll\'ar for
persuading us to rewrite the introduction.

\section{Birational maps} \label{sec:sbm}

In this section we explain several constructions of birational maps.
First of all we recall the definition of link of the Sarkisov program,
in a context which is sufficient for our purposes.
We continue with an informal discussion leading to a link 
$X_4 \dasharrow Y_{3,4}$; this is how we first discovered it. Then we 
give a much shorter, more efficient construction, which generalizes to many
other cases. This corresponds to ``Type I'' of \cite{CPR} and was explained 
to us by M. Reid. Finally, if $P\in L\subset X$ is a line on $X$, we construct
a link $X\dasharrow X$ centered on $L$; this presents some topical intricacies
which we treat quickly, because they are very similar to ``Type II'' 
of \cite{CPR}, 4.11 and 7.3.


\subsection{Links}

 \begin{dfn} \label{dfn:slk}
A {\em Sarkisov link of Type~II} between two Fano 3-folds $X$ and $Y$
is a birational map $f\colon X\dasharrow Y$ that
factorizes as
 \[
 \renewcommand{\arraystretch}{1.5}
 \begin{matrix}
 V&\dasharrow&V' \\
 \downarrow&&\downarrow \\
 X&\stackrel{f}{\dasharrow} & Y
 \end{matrix}
 \]
where
 \begin{enumerate}
 \renewcommand{\labelenumi}{(\alph{enumi})}
 \item $V\to X$ and $V'\to Y$ are extremal divisorial contractions in
the Mori category, and
 \item $V\dasharrow V^\prime$ is a composite of inverse flips, flops and flips
(in that order), and in particular, is an isomorphism in codimension $1$.
 \end{enumerate}
Usually (always in this paper) the map $V\dasharrow V^\prime$ is
a flop (flip, inverse flip):
 \[
 \diagram
   &V\dlto\drto\rrdashed|>\tip& &V^{\prime}\dlto\drto&    \\
 X &                          &Z&                    &Y\enddiagram
 \]
In this case we say that $Z$ is the {\em midpoint} of the link
 \end{dfn}

\subsection{Constructions, first approach}

Fix a quartic 3-fold $X=X_4 \subset\P^4$ with a singular point $P \in X$, 
locally analytically equivalent to the origin in the hypersurface
 \[
\{xy+z^3+t^3=0\}\subset \C^3
 \]
Changing coordinates, we may then write the equation of $X$ as:
 \begin{equation} \label{eq:1}
F=x_0^2x_1x_2+x_0a_3+b_4=0
 \end{equation}
where $a_3$, $b_4$ are homogeneous polynomials of degrees $3,4$ in the 
variables $x_1, x_2, x_3, x_4$.
In what follows we always assume that $X$ is general, in the sense that
 \begin{itemize}
 \label{conditions}
 \item[(a)] $X$ has only one singular point $P=(1:0:0:0:0)$,
 \item[(b)] $a_3(0,0,x_3,x_4)=b_4(0,0,x_3,x_4)=0$ only if $x_3=x_4=0$
 (this is a genericity condition involving the lines $P\in L \subset X$ 
 passing through $P$, see below).
 \end{itemize}

We begin with some heuristic considerations that lead to the construction of a 
birational map $X_4 \dasharrow Y_{3,4}$. Theorem~\ref{thm:lk1} states that 
this map is a link of the Sarkisov program. 

 \begin{lem} Let $X=\{xy+z^3+t^3=0\} \subset \C^4$ be a 3-fold germ. Let 
 $U\to \C^4$ be the weighted blowup with weights $(2,1,1,1)$, $V\subset U$ 
 the inverse image of $X$ and $E \subset V$ the exceptional divisor. Then

 \begin{enumerate}
 \renewcommand{\labelenumi}{(\arabic{enumi})}
 \item $E\subset V \to 0 \in X$ is an extremal divisorial contraction in 
 the Mori category with discrepancy $a_E(K_X)=1$,
 \item $E_{|E}=\O(-1)$, where $\O_E(1)$ denotes the tautological sheaf 
 under the obvious embedding 
 \[
 E=\{xy+z^3+t^3=0\} \subset \P(2,1^3).
 \]
 In particular, $E^3=3/2$. \qed
 \end{enumerate}
 \end{lem}

Later in Section~\ref{sec:dco}, we prove that any divisorial contraction 
$E\subset V \to P \in X$, contracting the exceptional 
divisor $E$ to $P$, is isomorphic to a weighted blowup with weights either 
$(2,1,1,1)$ or $(1,2,1,1)$.

Let $X=X_4$ be a quartic 3-fold as in equation~\ref{eq:1}, and let 
$E \subset V \to P\in X$ be the weighted blowup assigning weights 
$2,1,1,1$ to the coordinates $x_1, x_2, x_3, x_4$. We now play a 2-ray game
(see \cite{Co2}, pp. 269--272) starting with the configuration $V \to X$. 
We do this to determine if a link of the Sarkisov program exists originating
from this blowup. The heuristic calculations based on the Hilbert function 
are also explained in \cite{CPR}, 7.2. Denoting $A=-K_X$ and $B=-K_V$, we
have   
 \[
 B^3=(A-E)^3=A^3-E^3=4-\frac{3}{2}=\frac{5}{2}
 \]
Note that $N^1V=H^2(V, \mathbb{R})=\mathbb{R}^2$ is 2-dimensional, so 
$\NEbar V$ has two (pseudo)extremal rays; we denote the one, corresponding
to the curves contracted by $V \to X$, $R_{\text{old}}$, the other 
$R_{\text{new}}$. To perform the 2-ray game, the first step is to locate 
$R_{\text{new}}$ and determine its nature. Since $B^3>0$, we guess that $B$ 
is nef and that $R_{\text{new}}$ is a flop (in fact, with a little 
experience, this is easy to see: the functions $x_1,x_2,x_3,x_4$ all vanish 
on $E$ hence define sections of $B$, their common base locus is the 
union of the (proper transforms of the) 24 lines 
$P \in L_i \subset X$ given by $x_1x_2=a_3=b_4=0$. 
These divide in two groups of 12; those lying on $x_1=0$ have 
$B \cdot L_i =0$ and support the new ray, while those on $x_2=0$ have 
$B\cdot L_i=1/2$ and are thus not extremal), then $V$ is a weak Fano 3-fold 
(that is, $-K_V$ is nef and big) with anticanonical model 
 \[
 Z= \Proj \oplus_{n\geq 0} H^0(V,nB)
 \]

Quite generally if $Z$ is a Fano 3-fold with virtual singularities \linebreak 
$1/r_i (a_i, -a_i, 1)$, we
can write
\begin{eqnarray*}
-K_Z^3         &=&2g-2 +\sum \frac{a_i(r_i-a_i)}{r_i}   \\
h^0(Z, -K_Z)   &=&g+2
\end{eqnarray*}
(this uses the Riemann-Roch formula of Fletcher and Reid from \cite{YPG}; it 
is natural, by analogy with the classical case, to say that $g$ is the 
\textit{genus} of $Z$). In fact these numerical data determine the 
whole Hilbert function of $Z$
 \[
h^0(Z, -nK_Z)= -\frac{1}{12}n(n+1)(2n+1)K_Z^3+(2n+1)-\sum_{Q\in \Sing Z} 
\ell_Q (n+1)
 \] 
where the sum is taken over the virtual singularities $Q \in Z$ and the local 
contribution is given by the formula
 \[
\ell_Q (n)=\sum_{k=1}^{n-1}\frac{\overline{ka}(r-\overline{ka})}{2r}
 \]
(see \cite{YPG} for details, explanations and examples).

Now in our case $Z$ must be a Fano 3-fold with a singularity 
$1/2(1,1,1)$ and genus 2. This determines the Hilbert function of 
$Z$ uniquely, and we could try to use it to determine
$Z$, as explained in \cite{IF} or, more extensively, in \cite{k3}. However, 
it is easy to look first if we can spot  $Z$ in the list \cite{IF}
of (weighted) hypersurfaces or codimension 2 complete intersections. 
We easily find that $Z=Z_5 \subset \P(1,1,1,1,2)$ has the correct 
numerical invariants.
Let us try and construct a birational map $X\dasharrow \P(1,1,1,1,2)$ 
with image a variety of degree 5. We look for a section
$y \in H^0(V, 2B)$ where $B=A-E$; it is easy to see that
 \[
y=x_0x_1
 \]
gives such a section, indeed by construction $x_1$ vanishes twice along $E$. 
Finally
 \[
x_1F=y^2x_2+ya_3+x_1b_4=0
 \]
gives the equation of $Z$.
Next, the 2-ray game instructs us to flop the curves contracted by the morphism
$V \to Z$, that is, the proper transforms of the 12 lines 
$P \in L_i \subset X$ given by $x_1=a_3=b_4=0$.
Let $t\colon V \dasharrow V^\prime$ be the flop of these lines. 
Again $\NEbar V^\prime =R^\prime_{\text{old}}+R^\prime_{\text{new}}$ is a 
2-dimensional cone, with $R^\prime_{\text{old}}$ corresponding to the curves 
just flopped. As before,
we now wish to locate and determine the structure of $R^\prime_{\text{new}}$. 
Before we do this explicitly, we want to make a general remark. First, 
$V^\prime$ is uniruled hence $K_{V^\prime}$ is not nef. 
We know $K_{V^\prime}\cdot R^\prime_{\text{old}}=0$, therefore
$K_{V^\prime} \cdot R^\prime_{\text{new}}<0$. This means that the 2-ray game 
from now on is an ordinary minimal model program for $V^\prime$; in 
particular the existence of this minimal model program guarantees that the 
2-ray game ends in a link of the Sarkisov program. 

We can hope that the contraction of $R^\prime_{\text{new}}$  
is a divisorial contraction $S^\prime \subset V^\prime \to Y$, landing in a 
new Fano 3-fold $Y$ and completing a link $X\dasharrow Y$ of the Sarkisov 
program. To take this further, let us look for the exceptional surface 
$S^\prime$, or rather its transform $S\subset X$. We locate $S$ in $X$ as 
a \textit{special surface} \cite{CPR}, i.e. 
$S= \{f=0\}\cap X \subset X$ where $f\in k[x_0, ...,x_4]$ is 
a homogeneous function with highest slope
 \[
 \mu_P f =\frac{\mult_P f}{\deg f}
 \]
We already own $f=x_1$ with slope $\mu_P x_1=2$ and indeed it is easy to 
check that the proper transform $S^\prime \subset V^\prime$ of 
$X\cap \{x_1=0\}$ is a $\P^2$ with normal bundle $\O(-2)$. Putting
$A^\prime = -K_{Y}$, we have that $B^\prime = 
A^\prime -(1/2)S^\prime$ hence
 \[
 A^{\prime 3}=B^{\prime 3}+\frac{4}{8}=\frac{5}{2}+\frac{1}{2}=3
 \]
As before, this determines the Hilbert function of $Y$, a Fano 3-fold of 
genus $g=2$ and singularities $2 \times 1/2(1,1,1)$. We find $Y$ in the list 
of Fano 3-fold codimension 2 weighted complete intersections
 \[ 
 Y= Y_{3,4} \subset \P(1^4,2^2)
 \]
All that remains to do now is to exhibit the map 
$Z \dasharrow Y \subset \P(1^4,2^2)$ and the equations of $Y$. 
To do this observe that the strict transform of the cubic
$(x_0x_1x_2+a_3=0)_{|X}$ is in $|3B^\prime|$. The rational map is defined as 
follows
 \[
 Z \ni Q \mapsto \bigl( x_1(Q),..., x_4(Q), y_1(Q), y_2(Q) \bigr) 
 \in \P(1^4,2^2)
 \] 
where $y_1=y =x_0x_1$, and 
 \[
 y_2=-\frac{x_0x_1x_2+a_3}{x_1}
 \] 
Finally, it is easy to verify (exercise!) that $Z$ maps birationally onto its 
image $Y$, with equations
 \[
 \left\{\begin{array}{ll}
 y_1y_2+b_4(x_1,x_2,x_3,x_4) & = 0 \\ 
 y_1x_2+y_2x_1+a_3( x_1,x_2,x_3,x_4) & = 0
 \end{array}
 \right.
 \]
We have proved:
 \begin{thm} \label{thm:lk1} The diagram
 
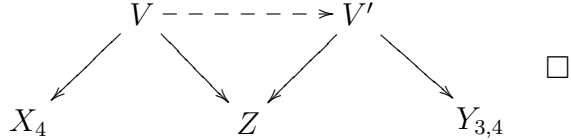
\[
 \diagram
    &V\dlto\drto\rrdashed|>\tip& &V^{\prime}\dlto\drto & \\
 X_4&                          &Z&               &Y_{3,4}\enddiagram
 \qed
 \]
is a Sarkisov link of Type~II
 \end{thm}

 \begin{rem} Note that there are really two links $X_4 \dasharrow
Y_{3,4}$, corresponding to the two weighted blowups of the singularity
$xy+z^3+t^3$, with weights $(2,1,1,1)$ and $(1,2,1,1)$. Note that both 
links go to the {\em same} variety $Y_{3,4}$
 \end{rem}

 \begin{exe} \label{exe:evil}  Construct the following links (recall
the Introduction, 1.5).
 \begin{enumerate}
 \renewcommand{\labelenumi}{(\alph{enumi})}
 \item Take $X=X_7 \subset \P(1,1,1,2,3)$ with a singular point at 
$P=(0,0,1,0,0)$ of the form $y^2+z^2+x_1^6+x_2^6$, then $X_7 \dasharrow 
Y_{6,7} \subset \P(1,1,2,3,3,4)$ with
midpoint $Z_9\subset \P(1,1,2,3,3)$.
 \item Take $X=X_{15} \subset \P(1,1,2,5,7)$ with a singular point at 
$P=(0,1,0,0,0)$
of the form $u^2+z^2y+y^7+x^{14}$, then 
$X_{15}\dasharrow Y_{14,15} \subset \P(1,2,5,6,7,9)$ with
midpoint $Z_{20}\subset \P(1,2,5,6,7)$.
 \end{enumerate}
 \end{exe}

\subsection{Constructions, second approach}

We briefly discuss a much more concise description of the link
$X_4 \dasharrow Y_{3,4}$ which also points out to a very large number
of similar examples. Everything here was suggested by Miles Reid.

The idea is to describe the link starting from the midpoint $Z$,
rather than either of the ends $X$, $Y$. In our example $Z$ is a 
Fano hypersurface in weighted projective space; it is special because
it is not $\Q$-factorial, this is just an expression of the fact that $Z$
is not $\Q$-factorial, that is, it is not in the Mori category.

Start with a Fano 3-fold hypersurface $Z=Z_d\subset \P(1,a_1,a_2,a_3,a_4)$
containing a surface
 \[
 \xi = \eta =0
 \]
Here $\xi, \eta$ are homogeneous functions of the coordinates of degrees
$\deg \xi < \deg \eta$. Usually one takes $\xi, \eta$ to be two 
coordinate functions, but not always. The equation of $Z$ can be written as
 \[
 F=a\eta-b\xi=0
 \]
Assume that $Z$ is quasismooth outside the ``nodes'' $a=b=\xi=\eta=0$, which
is often the case. Then we obtain two small partial resolution of $Z$, 
by considering the two ratios $y=\eta /\xi=b/a$ or 
$z=b/\eta=a/\xi$. 

For example, if $\xi=x_{i_1}$ and $\eta=x_{i_2}$ with $i_1<i_2$, then
the first ratio gives the hypersurface
 \[
 X=  \{ya(...x_{i_1},...x_{i_1}y,...) =b(...x_{i_1},...x_{i_1}y,...)\}
 \]
while the second ratio gives, in general, the codimension 2 complete
intersection
 \[
 Y:
 \begin{cases}
 z\eta = b \\
 z\xi = a
 \end{cases}
 \]
This construction explains all the quadratic involutions of
\cite{CPR}, our $X_4\dasharrow Y_{3,4}$ with midpoint 
$Z_5= \{a_3 y+b_4x_1=0\}$, and many more links involving
complete intersections $Y_{d_1,d_2}\subset \P(a_0, ..., a_5)$.

 \begin{exe} \label{exe:easy} Study the following cases
 \begin{enumerate}
 \renewcommand{\labelenumi}{(\alph{enumi})}
 \item $Z$ is a quartic 3-fold and $\xi=f_2$, $\eta=g_2$ are two quadrics,
 \item $Z$ is a quartic 3-fold and $\xi=x_1$, $\eta=g_2$,
 \item the examples of Exercise~\ref{exe:evil}
 \end{enumerate}
 \end{exe}

\subsection{Links centered on lines}

We next construct involutions \linebreak
$\tau_L \colon X \dasharrow X$ and show that they also are 
Sarkisov links of Type~II. These are very similar to the 
involutions in \cite{CPR} pp. 198 foll. ``Elliptic involutions''.

Let $L \subset X$ be a line passing through the singular point $P\in X$. 
Choosing coordinates so that $L=(x_2=x_3=x_4=0)$,
the equation of $X$ can be written as
 \[
F=x_2x_0^2x_1+a_1x_0x_1^2 +b_1x_1^3 +a_2x_0x_1
+b_2 x_1^2  +a_3 x_0 +b_3 x_1+b_4=0,
 \]
where $a_i$, $b_i$ are homogeneous polynomial of degree 
$i$ in $\C[x_2,x_3,x_4]$. It is easy to understand how a birational map 
$X \dasharrow X$ arises in this context: the generic fibre of the projection 
from $L$ is an elliptic curve with a section, corresponding
to the singular point $P \in X$. The birational map $X \dasharrow X$
is the reflection given by the group law on the elliptic curve. 
Our aim is to show that this birational selfmap is a link of the Sarkisov 
program. We follow closely the treatment of elliptic involutions in 
\cite{CPR}.

We eliminate both variables
$x_0$ and $x_1$ at once, replacing them by more complicated terms
 \[
 y=x_2x_0^2+a_1x_0x_1+b_1x_1^2+\cdots\quad\text{and}\quad 
 z=x_2x_0y+\cdots.
 \]
These are designed to be plurianticanonical on $V$, where 
$E\subset V \to L \subset X$ is the (unique) extremal divisorial contraction
which blows up the generic point of $L$. In other words, $y, z$ vanish 
enough times on the exceptional divisor $E$ of $V\to X$, and it turns out 
that, together with the other coordinates $x_2,x_3,x_4$, they generate
the anticanonical ring of $V$, and satisfy a relation of the form
 \begin{equation*}
 z^2+Azy+Bz=x_2y^3+Cy^2+Dy+E
 \end{equation*}
with $A,B,C,D,E\in k[x_2,x_3,x_4]$. This equation defines the
midpoint of the link, which is a (weak) Fano hypersurface
$Z_{10}\subset\P(1,1,1,3,5)$ having a biregular involution $i_Z$
coming from interchanging the two roots of the quadratic equation.

The form of the equation makes clear that the argument depends at some
level on the fact that the fibers of the rational map to $\P^2$
given by $x_2,x_3,x_4$ are birationally elliptic curves with a
section.

Define
 \[
 y= x_2x_0^2+a_1x_0x_1 +b_1x_1^2 +a_2x_0
 +b_2 x_1  +b_3
 \]
so that
 \begin{equation}\label{eq:2}
 F=yx_1+a_3x_0+b_4
 \end{equation}
From the last equation, it is clear that the divisor of zeros of $y$
on $V$ is $\geq 3E$ or, equivalently, $\mult_L y_{|X} = 3$. Next comes
the tricky step. Multiply F by $x_2x_0$, substitute for $x_2x_0^2$
in terms of $y$:
 \begin{align*}
 x_2x_0F&=x_2x_0x_1y+a_3x_2x_0^2+b_4x_2x_0\\
        &=x_2x_0x_1y+a_3(y-a_1x_0x_1-b_1x_1^2-\cdots) + b_4x_2x_0
 \end{align*}
Collecting the terms divisible by $x_1$ we then get
 \begin{multline*}
 x_2x_0F=x_1(x_2x_0y-a_1a_3x_0-a_3b_1x_1-a_3b_2)\\
 +a_3y-a_2a_3x_0-a_3b_3+b_4x_2x_0
 \end{multline*}
Set now
 \begin{equation}\label{eq:3}
 z=x_2x_0y-a_1a_3x_0-a_3b_1x_1-a_3b_2
 \end{equation}
so that
 \begin{equation} \label{eq:4}
 x_2x_0F=z x_1+ (x_2b_4-a_2a_3)x_0 +a_3(y-b_3)
 \end{equation}
Again the last equation makes it manifest that the divisor of zeros of $z$
on $V$ is $\geq 5E$ or, equivalently, $\mult_L z_{|X} = 5$. 
 
In order to eliminate $x_0,x_1$ in favor of $y,z$, note that we can view
(\ref{eq:2}--\ref{eq:4}) as inhomogeneous linear relations in $x_0,x_1$
with coefficients in $k[x_2,x_3,x_4,y,z]$:
 \[
 \renewcommand{\arraystretch}{1}
 \begin{array}{rrcccccccc}
 (\ref{eq:2})&F=      &yx_1 &+& a_3   x_0        &+&b_4       &=&0,\\
 (\ref{eq:3})&x_2x_0F=&z x_1&+&(x_2b_4-a_2a_3)x_0&+&a_3(y-b_3)&=&0,\\
 (\ref{eq:4})& \text{definition of $z$:}
&a_3b_1x_1&+&(a_1a_3-x_2y)x_0&+&z+a_3b_2                   &=&0.
 \end{array}
 \]
The equation relating $y$ and $z$ can be easily expressed
in the following determinantal \textit{format}
 \[
 \frac{1}{a_3}\det
 \begin{pmatrix} a_3 & y & b_4 \\
 x_2b_4-a_2a_3 & z & a_3(y-b_3)\\
 a_1a_3-x_2y   & a_3b_1 & z+a_3b_2 
 \end{pmatrix}
 =0.
 \]
The equation of $Z$ is quadratic in the 
last variable $z$, so that $Z$ is a 2-to-1 cover of
$\P(1^3, 3)$, which gives a (biregular!) involution of $i_L\colon
Z \to Z$ by interchanging
the sheets. The involution $\tau_L: X \dasharrow X$, as in \cite{CPR}, is the 
composite
 \[
 \renewcommand{\arraystretch}{1.3}
 \renewcommand{\arraycolsep}{0pt}
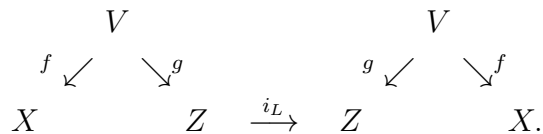
 \begin{matrix}
 && \ V \ &&&& \ V \ \\
 &{}^f\swarrow && \searrow^g &&{}^g\swarrow && \searrow^f \\
 X &&&&Z\quad \stackrel{i_L}\longrightarrow \quad Z&&&& X.
 \end{matrix}
 \]
 \begin{thm} \label{thm:lk2} The map 
$\tau_L \colon X \dasharrow X$ just constructed is a link of the Sarkisov 
program.
 \end{thm}
 \begin{proof} We have to show that $V \to Z$ contracts a finite
number of curves. The verification is tedious, but similar to 
\cite{CPR}, pp. 200--201.
 \end{proof}

 \begin{exe} Let $X$ be a quartic 3-fold, and $L\subset X$ a line on $X$.
If $X$ has one or more singular points along $L$, we get an elliptic
fibration and a link centered on $L$. It is 
very tricky, and amusing, to determine the structure of this link.
As an exercise, prove that:
 \begin{enumerate}
 \renewcommand{\labelenumi}{(\alph{enumi})}
\item if $X$ has one node on $L$, then $X \dasharrow X$ with midpoint
$Z_{12} \subset \P(1,1,1,4,6)$,
\item if $X$ has two nodes on $L$, then $X \dasharrow X$ with midpoint
$Z_{8} \subset \P(1,1,1,2,4)$
 \end{enumerate}
 \end{exe}

\section{Divisorial contractions} \label{sec:dco}

We survey the known results on the classification of 3-fold divisorial
contractions. Our main goal is to classify divisorial contractions 
contracting a divisor to the singularity $xy+z^3+t^3=0$. 

 \begin{dfn} \label{dfn:cont} Let $P\in X$ be the germ of a 3-fold
terminal singularity. A {\em divisorial contraction} is a proper birational
morphism $f\colon Y\to X$ such that
 \begin{enumerate}
 \renewcommand{\labelenumi}{(\arabic{enumi})}
 \item $Y$ has terminal singularities,
 \item the exceptional set of $f$ is an irreducible divisor $E\subset Y$,
 \item $-K_Y$ is relatively ample for $f$.
 \end{enumerate}

An {\em extremal} divisorial contraction $f\colon Y\to X$ is an extremal
divisorial contraction in the Mori category. In other words, $Y$ has
$\Q$-factorial terminal singularities, $f$ is the contraction of an extremal
ray $R$ of $\NEbar Y$ satisfying \hbox{$K_Y\cdot R<0$}, and the
exceptional set $\Exc f=E\subset Y$ is a divisor in $Y$. Its image
$\Ga=f(E)$ is a closed point or a curve of $X$, and we usually write
$f\colon(E\subset Y)\to(\Ga\subset X)$. Here $X$ is not necessarily a germ,
but $Y\to X$ is a divisorial contraction in the above sense above the germ
around any point $P\in\Ga$. Viewed from $X$, we also say that $f$ is an {\em
extremal extraction}, or that it {\em extracts} the valuation $v=v_E$ of
$k(X)$ from its {\em center} $\Ga=C(X,v_E)\subset X$.
 \end{dfn}

The classification of 3-fold divisorial contractions is now known in several 
important special cases:

 \begin{thm}[Kawamata \cite{Ka}] \label{thm:kaw} Let
 \[
 P\in X\cong \frac{1}{r}(1,a,r-a) \quad
 \text{\em (with $r\ge2$ and $a$ coprime to $r$)}
 \]
 be the germ of a $3$-fold terminal quotient singularity, and
 \[
 f\colon(E\subset Y)\to(\Ga\subset X)
 \]
 a divisorial contraction such that\/ $P\in\Ga$. Then\/ $\Ga=P$ and\/ $f$ is
 the weighted blowup with weights $(1,a,r-a)$. \qed
 \end{thm}

 \begin{cor} \label{cor:Kaw} Suppose that $X$ is a\/ $3$-fold with only
 terminal quotient singularities. If a curve\/ $\Ga\subset X$ is the center
 of a divisorial extraction\/ $f\colon(E\subset Y)\to(\Ga\subset X)$ then\/
 $\Ga\subset\NonSing X$ (and\/ $f$ is the blowup of\/ $I_\Ga$ over the
 generic point of $\Ga$).
 \end{cor}
 \begin{proof}
For if $\Ga$ passed through a terminal quotient point $P$,
Theorem~\ref{thm:kaw} would imply that $\Ga=P$, a contradiction. 
 \end{proof}
The next Corollary is a very useful characterization of terminal
singularities of pairs in terms of multiplicity.
 \begin{cor} \label{cor:multk} Let
 \[
 P\in X\cong \frac{1}{r}(1,a,r-a) \quad
 \text{\em (with $r\ge2$ and $a$ coprime to $r$)}
 \]
 be the germ of a $3$-fold terminal quotient singularity, and let $\H$
 be a linear system (not necessarily mobile) on $X$. Let 
 \[
 f\colon(E\subset Y)\to(P \in X)
 \]
 be the blow up with weights $(1,a,r-a)$, and $\delta=\mult_E \H$.
 Then 
 \[
 K+\frac{1}{n} \H
 \]
 is terminal, if and only if $\delta < n/r$.      
 \end{cor}
 \begin{proof} If $(X,(1/n) \H)$ is not terminal, by \cite{Co1}~\S2, there
   is a divisorial contraction $f\colon(E\subset Y)\to(P \in X)$,
   extracting a divisor $E$ for which $\mult_E \H \geq na_E$. By
   \ref{thm:kaw}, $f$ is the weighted blow up we are speaking of, so
   we get that $\delta \geq n /r$. Vice-versa if $(X, \H)$ is terminal,
   it is part of the definition that $\delta < n/r$.
 \end{proof} 

Another known case is in \cite{Co2}:

 \begin{thm}\label{thm:corti} Let $P \in X$ be a 3-fold germ analytically 
 isomorphic to an ordinary node
 \[
 xy+zt=0
 \]
 and $f \colon (E\subset Y) \to (P \in X)$ a divisorial contraction;
 assume in addition that  $f(E)=P$. Then 
 $f$ is the blow up of the maximal ideal at $P$. \qed 
 \end{thm}
The following Corollary is similar to \ref{cor:multk} but slightly
weaker, essentially because curve maximal centers passing through an
ordinary node do exist. 
 \begin{cor} \label{cor:multc} Let $P \in X$ be a 3-fold germ analytically 
 isomorphic to an ordinary node
 \[
 xy+zt=0
 \]
 and $\H$ a linear system (not necessarily mobile) on $X$. Let 
 \[
 f\colon(E\subset Y)\to(P \in X)
 \]
 be the blow up of the maximal ideal at $P$, and $d=\mult_E
 \H$. Assume that there exists a valuation $F$, with center $C_F X=P$,
 and $\mult_F \H >na_F$. Then $d>n$.
 \end{cor}
 \begin{proof} It is awkward to try to prove this by the same method
 as Corollary~\ref{cor:multk}, because of possible curve maximal
 centers. Fortunately, you can check that the proof of
 Theorem~\ref{thm:corti} in \cite{Co2}~pg.~282, proves the statement.
 \end{proof}

 Before stating our main result in this section, we mention the following
 very nice result of Kawakita \cite{Kw} 

 \begin{thm} Let $(E\subset Y) \overset{f}{\to} P \in X$ be a 3-fold 
divisorial contraction. Assume $P \in X$ is a nonsingular point and $f(E)=P$. 
In suitable analytic coordinates on $X$, $f$ is a weighted blow 
up; the weights are of the form $(1,a,b)$ with $\hcf (a,b)=1$. \qed
 \end{thm}

Our main result is the following:

 \begin{thm} \label{thm:dco} Let $x \in X$ be a 3-fold germ analytically 
isomorphic to
 \[
 xy+z^3+t^3=0
 \]
and $p \colon (E\subset Y) \to (x \in X)$ a divisorial contraction. Then 
$p$ is the weighted blow up with weights $(2,1,1,1)$ or $(1,2,1,1)$.
 \end{thm}

Before starting the proof, which occupies the rest of the section, we
like to make a few comments.

Our proof is a systematic integration of the ideas of \cite{Ka} and
\cite{Co2}, \S 3.2 and is not difficult in principle. However, it is rather 
complicated and it requires a great deal of notation. The main idea is to 
apply the connectedness theorem of Shokurov to many different morphisms. 
The reason all these morphisms exist is that the (analytic)
class group $\Cl (X,x)$ is large. We use this
manoeuvre many times and without much explanation; \cite{Co2}, Theorem 3.10
shows the technique at work in a much simpler situation.

It is certain that the method can be applied to other singularities,
for instance $xy +t^n+z^n$ and similar cases with relatively large local
class group. On the other hand, we have not been able to make the method
work when the class group is small; for instance, we have had little success 
so far with $xy+z^2+t^3$ (which is factorial). We have not attempted to 
determine how far the method can be pushed.

On the other hand a large part of the paper of Kawakita \cite{Kw}
holds for arbitrary cDV points $x \in X$. He was able recently to
extend his method to the classification of  divisorial
contractions to $xy+z^2+t^n$ \cite{Kw2}, and even general $cA_n$
singularities \cite{Kw3}. This includes our result as a special case.

\begin{proof} We begin by describing the general setup for the proof. 
\paragraph*{\textsc{General setup}}
Let $n$ be a sufficiently large and divisible positive integer;
fix a finite dimensional very ample linear system
 \[
 \H^Y \subset |-n K^Y|
 \]
Note that we write $K^Y$ to signify the canonical class of $Y$, as opposed
to the usual notation $K_Y$. We use the notation $K^Y$ throughout this proof.
Denote $\H =p(\H^Y)$ the image of $\H^Y$ in $X$, so that 
 \[
 K^Y+\frac{1}{n} \H^Y = p^\ast \Bigl(K^X + \frac{1}{n} \H \Bigr)
 \]
By construction
 \[
 \mult_E(\H)=n a_E(K_X)
 \]
while 
 \[
 \mult_\nu (\H)<n a_\nu (K_X)
 \]
for all valuations $\nu \not = E$.

\paragraph{\textsc{Description}} We summarize in the following
diagram, and explain below, our notation for the various spaces and
morphisms which we use in the course of the proof
 
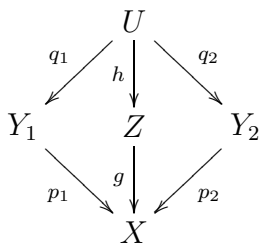
\[
 \diagram
              &U\dlto_{q_1}\dto_h\drto^{q_2}&                \\
 Y_1\drto_{p_1}&Z\dto_g                     &Y_2\dlto^{p_2}  \\
              &X                           &\enddiagram 
 \]
We now introduce in detail the various spaces and morphisms. Please draw your 
own picture, and do your own calculations to justify the 
description we give; otherwise you will not follow the proof.

\subparagraph{(a)} Denote $p_1 \colon E_1 \subset Y_1 \to x \in X$ the
weighted blow up with weights $(2,1,1,1)$ and exceptional divisor 
$E_1\subset Y_1$. Similarly denote $p_2 \colon E_2 \subset Y_2 \to x
\in X$ the weighted blow up 
with weights $(1,2,1,1)$. It is easy to check, for instance by
performing the blowing up explicitly, that $Y_i$ has a singular point 
$y_i \in Y_i$ of type $1/2(1,1,1)$ and is elsewhere nonsingular.

\subparagraph{(b)} Denote $g \colon Z \to X$ the blow up of the maximal
ideal at $x \in X$ with exceptional divisors $E_1$, $E_2$. The abuse
of notation means to suggest, for instance, that the rational map 
$Z\dasharrow Y_1$ (not a morphism!) is an isomorphism at the generic
point of $E_1$, and contracts $E_2$. Thus
the notation ``$E_1$'' denotes the ``same'' divisor in two different varieties
$Y_1$ and $Z$. We let the context decide which is meant; however, when we 
wish to be precise about the ambient variety, we write for instance $E_1^Z$ 
meaning the
divisor $E_1$ on the variety $Z$. We do this for other varieties and divisors, 
as well. This is justified since many of the quantities we are interested
in, such as discrepancies and multiplicities, depend only on the divisor, not
on the ambient variety.

\subparagraph{(c)} It is easy to check that $E_i^Z\cong \P^2$, and 
$E_1^Z$, $E_2^Z$ intersect in a line $B=E_1^Z \cap E_2^Z$. Also, $Z$
itself is nonsingular apart from three distinct ordinary nodes 
 \[
z_j \in B \subset E_1^Z+E_2^Z \subset Z , 
 \]
$j\in \{1,2,3\}$, each looking like
 \[ 
\text{origin}\in z\text{-axis}
\subset (t=0)  \subset (xy+zt=0) 
 \]

\subparagraph{(d)} Denote $h:U \to Z$ the blow up of the three $z_j$s; it has 
three exceptional divisors $F_j \subset U$, all isomorphic
to $\P^1 \times \P^1$ with normal bundles $N_{F_j} U \cong \O (-1,-1)$.
It is easy to check that $U$ is nonsingular. We also denote
$q_i \colon U \to Y_i$ the obvious morphisms and 
\[
f=g\circ h \colon U \to X
\]

\subparagraph{(e)} It is important to understand that the maps 
$Z \dasharrow Y_i$ are not morphisms. Indeed for instance we can
resolve the map $Z \dasharrow Y_1$ by a diagram 
$$\diagram
                &V_1\dlto \drto&   \\
Z\rrdashed|>\tip&              &Y_1\enddiagram$$
where $V_1 \to Z$ is a small resolution of all the $z_j\in Z$; 
$E_1^{V_1}\to E_1^Z$ blows up all the three points $z_j\in E_1^Z$, 
introducing exceptional curves $\Gamma^1_j = C_{F_j} (V_1)$, while
$E_2^{V_1}\to E_2^Z$ is an isomorphism. The morphism 
$V_1 \to Y_1$ contracts $E_2^{V_1} \cong \P^2$, with normal bundle
$N_{E_2}V_1 =\O(-2)$, to a singular point $1/2(1,1,1)$. The images
of the $\Gamma^1_j$ are three lines $C^1_j=C_{F_j} Y_1$ passing
through $y_1\in Y_1$. We use these lines later in the proof. Similar remarks
and notation apply to the map $Z \dasharrow Y_2$.

\paragraph{\textsc{Main division into cases}}
Now we start with our given divisorial contraction $E \subset Y \to x \in X$
and we want to show that $E= E_1$ or $E=E_2$; assuming the contrary we will
derive a contradiction. The proof divides out in cases, depending on the 
position of the center $C_E(Z)$ of $E$ on $Z$:
 \[
\text{{\sc Cases:}}
\begin{cases}
(1)    \quad &C_E(Z)\not \subset B=E_1^Z\cap E_2^Z            \\
(2.1)  \quad &C_E(Z)\subset B, \text{but} \; C_E(Z)\not = z_j \\
(2.2)  \quad &C_E(Z)=z_1 
\end{cases}
 \]
In all cases, we define $b_i$, $c_j$ by the formula
 \[
f^\ast \H = \H^U + \sum b_i E_i^U + \sum c_j F_j
 \]
Another way to say this is that $b_i =\mult_{E_i} \H$ and
$c_j= \mult_{F_j} \H$ are the multiplicities of $\H$ along
the $E_i$s and $F_j$s. The assumption that $E\not = E_1,E_2$ means that
$b_1, b_2<n$.

We treat each case separately.
\paragraph*{\textsc{Case 1}} Let $S\subset X$ be a generic surface
through $x \in X$. It is easy to compute
 \begin{align*}
f^\ast K^X =& K^U -\sum E_i^U -2 \sum F_j \\
f^\ast S   =& S^U +\sum E_i^U +  \sum F_j 
 \end{align*}
hence
 \begin{multline*}
 f^\ast \Bigl(K^X +S +\frac{1}{n}\H\Bigr)= K^U+S^U+\frac{1}{n}\H^U \\
 +\sum \frac{b_i}{n} E_i^U 
 +\sum \Bigl(\frac{c_j}{n}-1\Bigr)F_j
 \end{multline*}
which is another way to say, for instance, that 
$a_{E_i}(K^X +S +(1/n)\H)=b_i/n$, and similarly for the discrepancies of
the $F_j$.

Suppose now that $C_E Z \in E_1$, say. We apply Shokurov connectedness
theorem, \cite{Co2} \S 3.2 and 3.3 and especially Corollary~3.5,
to the morphism 
$f\colon U \to X$ and the divisor
 \[
K^U+S^U+ D^U=K^U+S^U+\frac{1}{n}\H^U  +\sum \frac{b_i}{n} E_i^U 
+\sum \Bigl(\frac{c_j}{n}-1\Bigr)F_j,
 \]
(where the formula defines $D^U$) and we conclude that there is a ``line'' 
$L_1 \subset E_1^U$ (by this we mean that $L_1$ maps to a line under
the morphism $E_1^U \to E_1^Z\cong \P^2$) such that
 \[
L_1 \subset LC \bigl(U, K^U+S^U+D^U\bigr)
 \]
One of the three $z_i$s, say $z_1$, does not lie on the line $L_1$.
It is easy to construct a contraction 
$W_1 \to X$ having $F_1$ as its unique exceptional divisor ($W_1$ has
canonical but not terminal singularities so this is not a divisorial
contraction in the Mori category). In fact there is a morphism 
$\varphi \colon U \to W_1$, resulting from a free linear system
$\Delta^U$ on $U$, which in turn is the proper transform of a linear
system $\Delta$ on $X$ such that 
\[
f^*\Delta=\Delta^U+\sum E_i+2F_1+\sum_{j=2,3}F_j
\]
Define now
\[
\delta=\max \Bigl\{ 0,\frac{1}{2}\Bigl(1-\frac{c_1}n \Bigr) \Bigr\}
\]
Then 
 \begin{multline*}
f^*\Bigl(K^X+S+\frac1n\H+\delta\Delta\Bigr)=K^U+S^U+\frac1n\H^U+\delta
\Delta^U+\\
\sum\Bigl(\frac{b_i}n+\delta\Bigr)E_i+p
F_1+\sum_{j=2,3}\Bigl(\frac{c_j}n+\delta-1\Bigr)F_j
 \end{multline*}
where $p=(c_1/n) -1 +2\delta \geq 0$. We now apply Shokurov
connectedness to the morphism $\varphi$ and the divisor
 \begin{multline*}K^U+S^U+D_1^U=K^U+S^U+\frac1n\H^U+\delta\Delta^U+\\
\sum \Bigl(\frac{b_i}n+\delta \Bigr)E_i+pF_1+\sum_{j=2,3}
\Bigl(\frac{c_j}n+\delta-1 \Bigr)F_j
 \end{multline*}
(where the formula defines $D_1^U$). It is important here to
understand that the divisor $F_1$ is not contracted by $\varphi 
\colon U \to W_1$; the reason we can apply
Shokurov connectedness is that $p\geq 0$. This is why we introduced
$\delta$ in the first place; it would have been tempting to run the
argument with $\delta=0$. Now $E_1^U$ is contracted by $\varphi$; 
the fibers of the rational map $\varphi_{|E_1^Z} \colon E_1^Z \cong
\P^2 \dasharrow W_1$ are the lines trough $z_1$. Furthermore
\[
L_1 \cup S^U\subset LC \bigl(U, K^U+S^U+D_1^U\bigr)
\]
Since $C_{F_1}(Y_1)$ is a curve, $a_{F_1}(K^{Y_1})=1$, and a simple
calculation then shows that $c_1=b_1+\mult_{F_1} \H^{Y_1}\geq b_1$,
therefore 
\[
\frac{b_1}n+\delta<1.
\]
This implies that $LC(U, K^U+S^U+D_1^U)$ is not connected in the neighborhood
of a general fiber of $\varphi_{|E_1^U} \colon E_1^U \to W_1$. The 
contradiction finishes the proof in Case 1.
\paragraph*{\textsc{Case 2}} Let us define $\delta_1$, $\delta_2$ by
the formulas
 \begin{align*}
q_1^\ast \H^{Y_1}&=\H^U + \delta_1 E_2^U +\text{(other)}\\
q_2^\ast \H^{Y_2}&=\H^U + \delta_2 E_1^U +\text{(other)}
 \end{align*}
where ``(other)'' means a combination of the $F_j$s. Another way to
define these numbers would have been to set
 \begin{align*}
\delta_1 = &\mult_{E_2} \H^{Y_1}\\ 
\delta_2 = &\mult_{E_1} \H^{Y_2} 
 \end{align*}
It is convenient (and, ultimately, straightforward) for us to calculate 
$\delta_1$ and $\delta_2$ in terms of $b_1$, $b_2$, in fact
we only need the

\paragraph*{\textsc{Claim}} 
 \[
\delta_1 + \delta_2 =\frac{b_1+b_2}{2}
 \]
To prove the claim, note
 \begin{align*}
q_1^\ast E_1 &= E_1^U + \frac{1}{2} E_2^U+\text{(other)} \\
q_2^\ast E_2 &= E_2^U + \frac{1}{2} E_1^U+\text{(other)} 
 \end{align*}
therefore
 \begin{multline*}
f^\ast \H =q_1^\ast p_1^\ast \H =q_1^\ast(\H^{Y_1}+b_1E_1)\\
=\H^U+b_1E_1^U+ \Big(\delta_1 +\frac{b_1}{2}\Bigr)E_2^U+\text{(other)}
 \end{multline*}
from which we conclude $b_2=\delta_1 + b_1/2$. Similarly,
$b_1=\delta_2 + b_2/2$ and the claim follows.

\paragraph*{\textsc{Subcase 2.1}} Assume that $C_E(Z)\subset B$ but is
not one of the $z_j$s. This implies that $C_E(Y_i)=y_i \in Y_i$ is the
unique singular point. The divisor
 \[
p_i^\ast \Bigl(K^X +\frac{1}{n}\H\Bigr)=
K^{Y_i}+\frac{1}{n} \H^{Y_i}-\Bigl(1-\frac{b_i}{n}\Bigr)E_i
 \]
is \textit{strictly canonical}, that is, canonical but not terminal, in a 
neighborhood of $y_i \in Y_i$ (the valuation $\nu$ corresponding to
the exceptional divisor of the divisorial contraction which we have been
studying all this time, has discrepancy $=0$). Because $b_i <n$, the divisor
 \[
K^{Y_i}+\frac{1}{n} \H^{Y_i}
 \]
is not canonical in a neighborhood of $y_i \in Y_i$. By 
Corollary~\ref{cor:multk}, both $\delta_1$, $\delta_2$ are $>n/2$ 
hence $(b_1+b_2)/2 > n$ and either $b_1$ or $b_2>n$, a contradiction
which concludes this case. 

\paragraph*{\textsc{Subcase 2.2}} Assume now that $C_E(Z)$ is one of
the $z_j$s, say $z_1$. The proof just given 
breaks down because typically in this case the center $C_E(Y_i)$ of $E$ 
on $Y_i$ is not the singular point $y_i \in Y_i$. Therefore we argue
directly on $Z$. It is important to be aware that the divisors $E_1^Z$
and $E_2^Z$ are not $\Q$-Cartier at $z_1$, but the sum $E_1^Z+E_2^Z$ is.
Consider the divisor
 \[
D^Z = \frac{1}{n} \H^Z +\frac{b_2-b_1}{n}E_2^Z
 \]
on $Z$. Note
 \begin{align*}
g^\ast \Bigl(K^X +\frac{1}{n} \H\Bigr) & = K^Z+\frac{1}{n}\H^Z
+\Bigl(\frac{b_1}{n}-1\Bigr)E_1^Z +\Bigl(\frac{b_2}{n}-1\Bigr)E_2^Z\\
 &=K^Z+D^Z+\Bigl(\frac{b_1}{n}-1\Bigr)(E_1^Z+E_2^Z)
 \end{align*}
hence $K^Z+D^Z$ is $\Q$-Cartier but not canonical at $C_EZ=z_1 \in Z$, hence
assuming as we may that $b_2\geq b_1$ so that $D^Z\geq 0$, and 
using Corollary~\ref{cor:multc}, we have
 \[
d=\mult_{F_1} D^Z >1
 \]
We now calculate $K^X+(1/n)\H$ in two different ways. On one hand
 \[
f^\ast \Bigl(K^X+\frac{1}{n}\H\Bigr)=K^U+\frac{1}{n}\H^U+
\Bigl(d+\frac{b_1}{n}-2\Bigr)F_1+\text{(other)}
 \]
while, on the other hand, using $f=p_iq_i$:
 \begin{multline*}
q_i^\ast p_i^\ast \Bigl(K^X+\frac{1}{n}\H\Bigr)
=q_i^\ast \Bigl(K^{Y_i}+\frac{1}{n}\H^{Y_i}+\frac{b_i}{n}E_i\Bigr)\\
=K^U+\frac{1}{n}\H^U+\frac{b_i}{n}E_i^U + \Bigl(\frac{\kappa_i}{n}
+\frac{b_i}{n}-2\Bigr)F_1+\text{(other)}
 \end{multline*}
where $\kappa_i=\mult_{F_1} \H^{Y_i}$; hence $\kappa_1/n=d>1$ gives
 \[
\kappa_1 > n
 \]
Similarly $\kappa_2/n+b_2/n-2 = d+b_1/n -2$ gives
 \[
\kappa_2+b_2-b_1 > n
 \]
Now recall the curves $y_i \in C^i_j=C_{F_j}Y_i \subset Y_i$;
by definition $\delta_1=\mult_{E_2} \H^{Y_1}$, and similarly for
$\delta_2$, while $\kappa_i =
\mult_{C^i_1}\H^{Y_i}$. It follows from Theorem \ref{thm:kaw} that
$\kappa_i/2\leq \delta_i$, hence combining inequalities we get
 \[
\frac{b_1+b_2}{2}=\delta_1+\delta_2\geq\frac{\kappa_1+\kappa_2}{2}
> n +\frac{b_1-b_2}{2}
 \]
and from this we conclude $b_2>n$, a contradiction.
\end{proof}

\section{Plan of proof of the main theorem} \label{sec:plan}

In this Section we give a more precise statement of the main theorem 
\ref{thm:main} and an outline of the proof. This is almost word by word
as \cite{CPR} \S 3. The proof is a formal consequence of the machinery and
definitions of the Sarkisov program, and the classification of
maximal singularities on $X_4$ and $Y_{3,4}$. We state the relevant results 
here, and prove them in sections \ref{sec:x4} and \ref{sec:x34}. Our precise 
statement is

 \begin{thm} \label{thm:main2} Let $X=X_4 \subset \P^4$ be a quartic 3-fold as 
above. In other words, $X$ has a singular point $P\in X$ of the form
$xy+z^3+t^3=0$, and is otherwise general. Let $V/T$ be an arbitrary
Mori fibre space and $\fie \colon X \dasharrow V$ a birational map. 
Then $\fie$ is a composition of the following 
birational maps:

 \begin{enumerate}
 \renewcommand{\labelenumi}{(\alph{enumi})}
 \item an involution $\tau_L \colon X \dasharrow X$, centered on a line
 $P\in L \subset X$ as in Theorem~\ref{thm:lk2}, 
 \item one of the two links $X\dasharrow Y_{3,4}$ as in Theorem~\ref{thm:lk1},

 \item The inverse of (b).

 \end{enumerate}
In particular, $V$ is isomorphic to either $X_4$, or $Y_{3,4}$, 
and $\mathcal{P}(X)=\{X_4,Y_{3,4}\}$.
 \end{thm}

Before giving an outline of the proof, we quickly recall some basic notions 
from the Sarkisov program. We refer to \cite{CPR} Chapter 3 for more details 
and discussion on this material.

 \begin{dfn}[canonical threshold] \label{dfn:c}
 $X$ is a variety, $\H$ a mobile linear system, and $\wave X\to X$ a
resolution with exceptional divisors $E_i$. As usual, we write
 \begin{align*}
 K_{\tilde X} &=K_X+\sum a_iE_i,\\
 \wave\H &=\H-\sum m_iE_i,
 \end{align*}
to define the {\em discrepancies} $a_i=a_{E_i}(K_X)$ of the exceptional
divisors $E_i$ and their {\em multi\-plicities} $m_i$ in the base locus of
$\H$. For $\la\in\Q$, we say that $K_X+\la\H$ is {\em canonical}\/ if all
$\la m_i\le a_i$, so that $K_{\tilde X}+\la\wave\H- (K_X+\la\H)$ is
effective ($\ge0$). Then we define the {\em canonical threshold} to be
 \begin{align*}
 c&=c(X,\H)=\max\bigl\{\la\bigm|\text{$K+\la\H$ is 
canonical}\bigr\}
 \\
 &=\min\nolimits_{E_i}\bigl\{a_i/m_i\bigr\}.
 \end{align*}
This is well defined, independently of the resolution $\wave X$. In all the
cases we're interested in, $K_X+(1/n)\H=0$ and $c<1/n$.
 \end{dfn}

 \begin{dfn}[maximal singularity] \label{dfn:ms}

Now suppose that \linebreak $K_X+(1/n)\H=0$ and $K_X+(1/n)\H$ is not
canonical, so that $c<1/n$. We make the following definitions:

 \begin{enumerate}
 \renewcommand{\labelenumi}{(\arabic{enumi})}
 \item a {\em weak maximal singularity} of $\H$ is a valuation $v_E$ of
$k(X)$ for which $m_E(\H)\ge na_E(K_X)$;

 \item a {\em maximal singularity} is an extremal extraction $Z\to X$
in the Mori category (see Definition~\ref{dfn:cont}) having exceptional
divisor $E$ with $c=a_E(K_X)/m_E(\H)$.

 \end{enumerate}
 \end{dfn}

In either case, the image of $E$ in $X$, or the center $C(X,v_E)$ of the
valuation $v_E$, is called the {\em center}\/ of the maximal singularity $E$.

 \begin{dfn}[degree of $\fie$] \label{dfn:deg}
 Suppose that $X$ is a Fano 3-fold with the property that $A=-K_X$ generates
the Weil divisor class group: $\WCl X=\Z\cdot A$ (this holds in our case).
Let $\fie\colon X\dasharrow V$ be a birational map to a
given Mori fibre space $V\to T$, and fix a very ample linear system $\H_V$
on $V$; write $\H=\H_X$ for the birational transform $\fie^{-1}_*(\H_V)$.

The {\em degree} of $\fie$, relative to the given $V$ and $\H_V$, is the
natural number $n=\deg\fie$ defined by $\H=nA$, or equivalently
$K_X+(1/n)\H=0$.
 \end{dfn}

 \begin{dfn}[untwisting] \label{dfn:un}
 Let $\fie\colon X\dasharrow V$ be a birational map as above, and $f\colon
X\dasharrow X'$ a Sarkisov link of Type~II. We say that $f$ {\em untwists}
$\fie$ if $\fie'=\fie\circ f^{-1} \colon X'\dasharrow V$ has degree smaller 
than $\fie$.
 \end{dfn}

 \begin{rem}
 The Sarkisov program factorizes an arbitrary birational map between Mori
fibers spaces as a chain of more general types of links, using a more
complicated inductive framework. See \cite{Co1}, Definition~3.4 for the
general definition of a {\em Sarkisov link} $f\colon X\dasharrow X'$, and
\cite{Co1}, Definition~5.1 for the {\em Sarkisov degree} of a birational map
$\fie\colon X\dasharrow Y$ between Mori fibre spaces. The above
Definitions are special cases that are sufficient
for our purposes in this paper. We can get away with this because we start
from our quartic $X=X_4$, and we only ever perform
untwistings that either return to $X$, or to the Fano 3-fold $Y_{3,4}$
 \end{rem}

\begin{lem} \label{lem:untw} Let $X$, $V/T$ be as before and
$\fie: X \dasharrow V$ a birational map. If $E\subset Z \to X$ is a maximal 
singularity, any Type~II link $X \dasharrow X'$ (as in 
Definition~\ref{dfn:slk}), starting with the extraction $Z\to X$,
untwists $\fie$.
\end{lem}

\begin{proof} \cite{CPR}, Lemma~4.2. 
\end{proof}

The following classification of maximal singularities on $X_4$ and $Y_{3,4}$ 
implies our main theorem
\ref{thm:main2}.

 \begin{thm} \label{thm:msx4} Let $X=X_4 \subset \P^4$ be a quartic 3-fold 
as in the assumption of Theorem~\ref{thm:main}, and $E$ a maximal 
singularity of $X$. Either:

 \begin{enumerate}
 \item $E\subset Z \to P\in X$ is one of the blow ups with weights 
$(2,1,1,1)$ or $(1,2,1,1)$ described above, or

 \item the center $C(E, X)=L$ is a line $P \in L \subset X$ and $E$ is 
generically the blow up of the ideal of $L$ in $X$. 

 \end{enumerate}   
 \end{thm}

 \begin{thm} \label{thm:msx34} Let $E$ be a maximal singularity on 
a general $Y_{3,4}\subset \P(1^4, 2^2)$, then 
$E \subset Z \to P \in Y$ is the standard blow up of either one of the two 
singular points of $Y$ on the line $x_1=x_2=x_3=x_4=0$ in $\P(1^4, 2^2)$.
 \end{thm}

This theorems summarize the conclusions of a whole series of calculations
carried out for \ref{thm:msx4} in Section~\ref{sec:x4} and for 
\ref{thm:msx34} in Section~\ref{sec:x34}. 

\begin{proof}[Proof that \ref{thm:msx4} and \ref{thm:msx34} imply 
Theorem~\ref{thm:main2}] 
This is standard, and is the same as the proof in \cite{CPR}. If $X$ is
Fano and $V\to T$ a Mori fibre space, a birational map $\fie\colon X\dasharrow
V$ is defined by a mobile linear system $\H$. By the
Norther-Fano-Iskovskikh inequalities
\cite{Co1}, Theorem~4.2, if $\fie$ is not an
isomorphism then $\H$ has a maximal center $P$ or $C$, hence a 
maximal singularity $E \subset Z\to P$ or $C$ by \cite{Co1},
Proposition~2.10. By Theorem~\ref{thm:msx4} and \ref{thm:msx34}, there is a 
birational map $i\colon X\dasharrow X^\prime$, where either $X'=X$ or
$X'=Y$, that is a Sarkisov link, and by 
Lemma~\ref{lem:untw} untwists the maximal center $P$ or $C$, so that 
$\fie\circ i$ has smaller degree. Thus after a
number of steps, either $X\cong V$ or $Y \cong V$. 
\end{proof}

\section{Excluding} \label{sec:excl}
Let $W$ be a center on a Fano 3-fold $X$; that is, $W=P\in X$
or $W=\Gamma \subset X$ is either a point or a curve on $X$. Eventually
in the next two Sections, we take $X=X_4 \subset \P^4$ our special singular 
quartic 3-fold, or $X=Y_{3,4}\subset \P(1^4, 2^2)$, but here we keep the 
discussion general.
We are concerned with the problem of ``excluding $W$'', that is, to prove
that $W$ is not a maximal center for \textit{any} linear system $\H$ on $X$.
In this Section we explain our general strategy for doing this. 

\subsection{Reduction to a surface problem} The first step is to reduce
to a surface question.

\subsubsection{The starting point}
 \begin{enumerate}
 \renewcommand{\labelenumi}{(\alph{enumi})}

 \item We assume by contradiction that $W$ is a maximal center: there is a 
mobile linear system $\H \subset |\O_X(n)|$ on $X$, and a valuation
$E$ with center $C_E X=W$ and $m_E \H > na_E K_X$. 

 \item We select a \textit{test linear system} $\T$ on $X$ with $W\subset
\bs \T$ contained in the base locus of $\T$. Often we take $\T =
|I_W (1)|$, but this does not always work. In the simplest cases, but not in
all cases, $W= \bs \T$. The choice of the test system is often delicate.
 \end{enumerate}
 
\subsubsection{The strategy} We work with a general member $S \in \T$;
the argument is slightly different according to whether the center $W$
is a curve or a point.
 \begin{description}
 \item[The center is a curve] 
 the assumption means that $\mult_W \H =m >n$, so we have
 \[
 \H_{|S}=\L+m'W+\text{(other)} \subset |\O_S(n)|
 \]
 where $\L$ is the mobile part of $\H_{|S}$. In general $m'\geq m$ but
 in most applications below $m'=m$. We concentrate on showing
 that the mobile system $\L$ can not exist. The idea of course is simple:
 a non empty linear subsystem in $\O_S(n)$ is unlikely to have a fixed
 part as large as $mW$, $m>n$. 

 \item[The center is a point] by construction
 \[
 K+S+\frac{1}{n}\H
 \]
 is not log canonical in a neighborhood of $W$. By Shokurov's inversion of 
 adjunction, see \cite{FA} 17.7
 \[
 K_S+\frac{1}{n}\H_{|S}
 \]
 is also not log canonical. Here the method works better if $\H_{|S}$ is 
 mobile
 but in general we have to allow $\H_{|S}=\L+F$ with nonempty fixed part $F$.
 We try to reach a contradiction by choosing general members
 $L_1, L_2$ in $\L$ and calculating the intersection number
 $L_1\cdot L_2$ on $S$. Theorem~\ref{thm_2ds} states that, if
 $K_S+(1/n)(\L+F)$ is not log canonical at $P$, then the local
 intersection number $(L_1\cdot L_2)_P$ at $P$ is large; the contradiction
 happens when it is too large. As before, the idea is very simple, even 
 crude: two curves in $\O_S(n)$ can not intersect in too many points.
 \end{description}

\subsection{Linear system on surfaces} The following theorem is very
 useful in the study of linear system on surfaces: 

\begin{thm} \label{thm_2ds}
Suppose that $P \in \Delta_1 + \Delta_2\subset S$ is the analytic germ
of a normal crossing curve on a nonsingular surface. Let $\L$ be a
mobile linear system on $S$ and denote $\L^2$ the local intersection
multiplicity $(L_1\cdot L_2)_P$ at $P$ of two general members $L_1$,
$L_2\in \L$. Fix rational numbers $a_1,a_2\ge0$ and suppose that
 \[
 K_S+(1-a_1)\Delta_1+(1-a_2)\Delta_2+\frac{1}{m} \L
 \]
is not log canonical for some $m>0$.
 \begin{enumerate}
 \item If either $a_1\le1$ or $a_2\le1$ then
 \[
 \L^2>4a_1a_2m^2.
 \]

 \item If both $a_i>1$ then
 \[
 \L^2>4(a_1+a_2-1)m^2.
 \]
 \end{enumerate}
\end{thm}

\begin{proof} In \cite{Co2}, Theorem~3.1.
\end{proof}

\subsection{The surface problem} We summarize the general set up 
for the surface problem in very broad terms.

\subsubsection{The set up} The initial set up is as follows.

 \begin{enumerate}
 \renewcommand{\labelenumi}{(\alph{enumi})}

 \item A polarized surface $S$, with polarizing (integral Weil, often Cartier) 
 divisor $A=\O_S(1)$. In practice, $S$ is often a K3 surface with DuVal 
 singularities.
 \item A configuration $\{\Gamma_i\}$ of curves $\Gamma_i \subset S$.  
 Usually $\Gamma_i$ is a -2 curve on the minimal resolution. The intersection 
 matrix $a_{ij}=\Gamma_i \cdot \Gamma_j$ is known in principle but it may be 
 impractical to calculate exactly, because there are many $\Gamma_i$ 
 or the geometry of the configuration maybe itself complicated or contain a 
 fair number of different degenerations. 
 Some geometric information is easily accessible, for instance the 
 $\Gamma_i$ are linearly independent in $H_2(S)$, any proper subset of 
 $\{\Gamma_i\}$ is contractible on $S$, $\NEbar S =\sum \Q_+[\Gamma_i]$ etc.  
 \item $A=\sum b_i \Gamma_i$ with $b_i$ small; often all $b_i=1$, sometimes 
 some $b_i=2$.
 \item We assume a nef $\Q$-divisor $L$ on $S$ given by a formula
 \[
 A=L + \sum \gamma_i \Gamma_i
 \]
 with the $\gamma_i \geq 0$. In the notation of the previous subsection, 
 this is $A=(1/n)\H_{|S}=(1/n)\L+ (1/n)F$, that is $L=(1/n)\L$
 and $\sum \gamma_i \Gamma_i=(1/n)F$.
 \end{enumerate}

\subsubsection{The goal} The aim is slightly different
 according to whether $W$ is a curve or a point.
 \begin{description}
 \item[When $W$ is a curve] the aim is to show that all $\gamma_i \leq 1$.
 \item[When $W$ is a point] assuming an inequality of the form
 \[
 L^2 > 4(2-\gamma_0)(1-\gamma_1),
 \]
 the aim is to show that $L$ cannot exist.
 \end{description}

\section{Centers on $X_4$} \label{sec:x4}

We fix $X=X_4 \subset \P^4$, with a singular point $P\in X$, satisfying
all the assumptions of Theorem~\ref{thm:main}. Our main goal in
this section is to prove the following:

 \begin{thm} \label{nocurves} A curve $\Gamma \subset X$, other than a line 
$P \in L \subset X$, can not be a maximal center. 
 \end{thm}

\begin{proof} Let $\Gamma$ be a curve and assume that $\Gamma$ is a
maximal center for $\H \subset |\O_X(n)|$. This implies
that $m=\mult_{\Gamma}\H> n$. In the proof, we reach a contradiction in
several steps:
\begin{description}
\item[Step 1] A raw argument shows that $\deg \Gamma \leq 3$.
\item[Step 2] $\Gamma$ can not be a space curve.
\item[Step 3] $\Gamma$ can not be a plane curve.  
\end{description}

\paragraph*{\textsc{Step 1}} Choosing general members $H_1$, $H_2$ of $\H$ and 
intersecting with a general hyperplane section $S$ we obtain
 \[
 4n^2=H_1\cdot H_2\cdot S>m^2\deg\Gamma.
 \]
This implies that $\deg \Gamma \leq 3$.

\paragraph*{\textsc{Step 2: space curves}} If $\Gamma$ is a space curve, then
by Step~1 it must be a rational normal curve of degree 3, contained in
a hyperplane $\Pi \cong \P^3 \subset \P^4$.
Let $S \in |I_{\Gamma,X}(2)|$ be a general quadric 
vanishing on $\Gamma$,$\L$ the mobile part of $\H_{|S}$; write 
 \[
 A = \O_S(1)=\frac{1}{n}\H_{|S}=L+\gamma\Gamma,
 \]
where $L=(1/n)\L$ is nef. Note that, because $I_\Gamma$ is cut out by
quadrics, 
\begin{displaymath}
  \mult_\Gamma \H = \mult_\Gamma \H_{|S} = n \gamma > n.
\end{displaymath}
We reach a contradiction by showing that $\gamma\leq 1$. For
simplicity we treat two separate cases, namely:
\begin{description}
\item[Case 2.1] $P \not \in \Gamma$,
\item[Case 2.2] $P \in \Gamma$.
\end{description}

\subparagraph*{\textsc{Case 2.1}} Here we assume that $P \not \in
\Gamma$. It follows that $S$ is nonsingular and $\Gamma^2=-5$ (all  
calculations of intersection numbers are performed on $S$).
Indeed it is easy to see that $S=S_{2,4} \subset \P^4$ 
is a nonsingular complete intersection of a quadric and
a quartic, therefore $K_S = \O_S(1)$. Then:
\begin{displaymath}
  -2 =\deg K_\Gamma = \Gamma \cdot (K_S + \Gamma ) = 3 + \Gamma^2
\end{displaymath}
shows that $\Gamma^2 = -5$. A simple calculation then gives:
\begin{displaymath}
  0\leq L^2=(A-\gamma \Gamma)^2=8-6\gamma -5 \gamma^2
\end{displaymath}
This shows that $\gamma\leq 4/5<1$ and finishes the proof in this
case. Note that we only need $\gamma \leq 1$; the additional room in
the argument, is what ultimately makes it possible to treat the next 
Case~2.2 essentially by the same method. 
 
\subparagraph*{\textsc{Case 2.2}} Now we assume that $P \in \Gamma$. 
Write as in Section~\ref{sec:dco} $\f \colon U\to X$ the resolution of 
singularities of $P \in X$, constructed in the proof of
Theorem~\ref{thm:dco}. Using that $I_{\Gamma ,\P^4}$ is 
generated by quadrics, and that a general quadric $Q\in |I_{\Gamma,
\P^4}(2)|$ is nonsingular, it is easy to see that the proper transform
$S^U$ must itself be nonsingular. 

Denote $\psi =\f_{|S^U}\colon S^U\to S$ and write 
\begin{equation*}
\psi^*\Gamma=\Gamma^U+\sum \gamma_i\Gamma_i+\sum g_j \Delta_j,
\end{equation*}
where $\Gamma_i=E^U_{i|S^U}$ and $\Delta_j=F_{j|S^U}$ are
$(-2)$-curves (here, following the notation of the proof of
Theorem~\ref{thm:dco}, $E_1^U$, $E_2^U$ and $F_1^U$, $F_2^U$, $F_3^U$
are the exceptional divisors of $U \to X$). 

There are now two subcases (up to relabelling the exceptional
divisors), depending on how the curve $\Gamma$ ``sits'' in the
singularity $P\in X$.
The crucial observation is that, because $\Gamma$ is a nonsingular
curve, $\Gamma^U$ intersect transversally a unique exceptional
divisor. The cases are as follows:
\begin{description}
\item[Subcase 2.2.1] The proper transform $\Gamma^U$ intersects
  $E_1^U$. In this case, $S^U$ meets $E_1^U$, $E_2^U$ and is disjoint
  from all the $F_j$s. Here $P\in S$ is an
  $A_2$-singularity, and $(\gamma_1,\gamma_2,g_1,g_2,g_3)=
  (2/3,1/3,0,0,0)$. 
\item[Subcase 2.2.2] The proper transform $\Gamma^U$ intersects
  $F_1^U$. In this case, $S^U$ meets $E_1^U$, $E_2^U$, $F_1^U$ and is 
  disjoint from $F_2$ and $F_3$. Here $P\in S$ is an
  $A_3$-singularity, and $(\gamma_1,\gamma_2,g_1,g_2,g_3)=
  (1/2,1/2,1,0,0)$. 
\end{description}

We claim that in both subcases
$$\Gamma^2\leq -4$$
Indeed it is easy to see, as in Case~2.1 ($P \not \in S$), that
$\Gamma^U\cdot \Gamma^U=-5$, and by the projection formula
\begin{displaymath}
  \Gamma \cdot \Gamma = \Gamma^U \cdot \psi^\ast \Gamma =
 \begin{cases}
  \Gamma^U\cdot \Gamma^U+ 2/3 =-5+2/3\\
  \Gamma^U\cdot \Gamma^U + 1 =-5+1
 \end{cases}
\end{displaymath}
in the two Subcases~2.2.1 and 2.2.2, respectively.
Finally
$$0\leq L^2=8-6\gamma + \Gamma^2 \gamma^2\leq 8-6\gamma -4\gamma^2$$
implies $\gamma<1$, a contradiction which concludes Step~2.

\paragraph*{\textsc{Step 3: plane curves}} Here we assume that
$\Gamma$ is a plane
curve of degree $d$ (by Step~1, $d\leq 3$), other than a line passing 
through $P$. Here too, as in Step~2, it is helpful and convenient to 
treat two cases, namely:
\begin{description}
\item[Case 3.1] $P \not \in \Gamma$ and $1\leq d \leq 3$.
\item[Case 3.2] $P \in \Gamma$ and $2\leq d \leq 3$.
\end{description}

\subparagraph*{\textsc{Case 3.1}} We first deal with the easy case 
$P\not\in \Gamma$ (following well known ``ancient'' methods of
Iskovskikh and Manin). Consider as usual a general element $S \in
|I_{\Gamma, X}(d)|$, denote $\L$ the mobile part of $\H_{|S}$; write 
 \[
 A= \O_S(1) =\frac{1}{n}\H_{|S}=L+\gamma\Gamma,
 \]
where $L=(1/n)\L$ is nef. We aim to show that $\gamma\leq 1$.
It is easy to see that $S=S_{d,4} \subset \P^4$ is a nonsingular
complete intersection of a quartic with
a hypersurface of degree $d\leq 3$, therefore $K_S =
\O_S(d-1)$. If $d\le 2$, then $p_a \Gamma =0$ and:
\begin{displaymath}
  -2 =\deg K_\Gamma = (\Gamma \cdot K_S + \Gamma ) = d(d-1) + \Gamma^2
\end{displaymath}
shows that $\Gamma^2 = -2-d(d-1)$. A simple calculation gives:
\begin{multline*}
  0\leq L^2=(A-\gamma \Gamma)^2=A^2-2 A \cdot \Gamma \gamma +\Gamma^2
   \gamma^2 \\
   = 4d-2d\gamma -d(d-1) \gamma^2-2\gamma^2
\end{multline*}
which implies that $\gamma \leq 1$. The proof is similar when $d=3$: 
$\Gamma^2=-6$, and then 
$0\leq L^2=(A-\gamma \Gamma)^2=A^2-2 A \cdot \Gamma \gamma +
\Gamma^2 \gamma^2 = 12-6\gamma -6 \gamma^2$ and again $\gamma \leq 1$.

These calculations finish Case~3.1 $P\not\in \Gamma$.

\subparagraph*{\textsc{Case 3.2}} From now on we assume that 
$P \in \Gamma$, $\Gamma$ 
not a line. In this case, restriction to a general element of the test
linear system $|I_\Gamma (d)|$ does not lead to a contradiction; it is 
necessary to use a different test system. 

Denote $\Pi \subset \P^4$ the plane spanned by $\Gamma$, let $S_1$,
$S_2$ be general hyperplane sections of $X$ containing $\Gamma$. 

We work with the ``test system'' $|S_1, S_2|$, even though
$\Gamma$ is usually only a component of its base locus $C=S_1\cap S_2 =
\bs |S_1, S_2|= X \cap \Pi$. We are assuming that $X$ is general,
hence in particular it is
terminal and $\mathbb{Q}$-factorial. This implies that $\Pi$ can not
be contained in $X$, and $C$ is a \emph{curve}.
Unfortunately, we have to divide the proof in several cases according
to what $C$ is. In the end each case is not very different or harder 
than any of the other cases, but we could not find a unified
presentation. The cases are as follows:
\begin{enumerate}
 \renewcommand{\labelenumi}{(\alph{enumi})}
 \item $C = \text{cubic} + \text{line}$, 
 \item $C = \text{conic} + 2 \, \text{lines}$,
 \item $C= 2 \, \text{conics}$,
 \item $C= \text{conic} + \text{double line}$,
 \item $C= \text{double conic}$. 
\end{enumerate}
We now treat Cases~(a)--(c); at the end we will show that Cases~(d)
and (e) do not happen (at least assuming, as we do, that $X$ is
general), in other words, we will show that $C=\Pi \cap X$ is always
\emph{reduced} when $X$ is general. Before treating each case
individually, we make some general comments and fix
the notation for the whole argument. 

Assuming for now that $C$ is reduced, we restrict to $S_1$ and write
\begin{eqnarray*}
A =(1/n)\H_{|S_1}& =L+\gamma\Gamma+\sum \gamma_i\Gamma_i \\
       S_{2|S_1} & =C=\Gamma+\sum \Gamma_i
\end{eqnarray*}
Our technique consists in selecting a ``most favorable'' component of
$C$, calculating an intersection number on $S_1$ using that $L$
is nef, and finally get that $\gamma\leq 1$.
When $C$ is reduced, it is clear that if $W$ is a component of $C$, 
$\mult_W \H=\mult_W\H_{|S_1}$; in particular $\mult_\Gamma \H
=\mult_\Gamma \H_{|S_1}=n\gamma$, and also $\mult_{\Gamma_i} \H =
n\gamma_i$. Because $\Gamma$ is a \emph{maximal} singularity, $\gamma 
\geq \gamma_i$, hence possibly after relabelling
components of $C$, we can assume that: 
 \[
 \gamma\geq \gamma_2\geq \gamma_1.
 \]
(ignore the term $\gamma_2$ if no curve $\Gamma_2$ is present).
Consider now the effective $\Q$-divisor \label{page:plane}
 \[
 (A-\gamma_1 S_2)_{|S_1}=L+(\gamma-\gamma_1)
 \Gamma+(\gamma_2-\gamma_1)\Gamma_2.
 \]
In Cases~(a) and (b), $\Gamma_1$ is a line and
 \[
 (1-\gamma_1)=(A-\gamma_1S_2)\cdot \Gamma_1 
 \geq (\gamma-\gamma_1)
 \Gamma\cdot \Gamma_1.
 \]
We now show that $\Gamma\cdot\Gamma_1\geq 1$ (on $S_1$); together with
the last displayed equation this implies that $\gamma \leq 1$ and
finishes the proof in Cases~(a) and (b). Note that this is intuitively
almost obvious: for example in Case~(a) $C$ is the plane union of a
cubic and a line, and we expect these to intersect in 3 points (when
we only need one!). The problem with saying this is, of course, that 
the set theoretic intersection $\Gamma \cap \Gamma_1$ can be all 
concentrated on the singular point $P\in X$. We now study this
situation more carefully. 

Note first that $S_1$ is nonsingular outside $P$. This follows easily
from the fact that the base locus $C$ of $|S_1,S_2|$ is a reduced
curve with only \emph{planar} singularities, and $X$ itself is
nonsingular outside of $P$ (this is all familiar and easy: if $f: Y
\to X$ is the blow up of $X$ along $C$, then $Y$ has isolated
singularities outside $f^{-1} (P)$).  

By our generality assumption~$2.2(b)$, and using the notation of the proof of 
Theorem~\ref{thm:dco}, we have that $\Gamma_1^Z \cap E_1^Z\cap E_2^Z=
\emptyset$. Also, $S^Z_{1|E_i}$ is nonsingular away from $E_1^Z\cap E_2^Z$.
Therefore \emph{either} the set theoretic intersection 
$\Gamma\cap \Gamma_1$ contains a nonsingular point of $S_1$,
\emph{or} $\Gamma^Z\cap\Gamma^Z_1$ contains a nonsingular point of 
$S^Z_1$. In both cases this point contributes with an integer value
$\geq 1$ to the intersection number $\Gamma \cdot \Gamma_1$, hence
our claim that $\Gamma\cdot\Gamma_1\geq 1$. This finishes the
proof in Cases~(a), (b).

In Case~(c), $\Gamma$ and $\Gamma_1$ are both conics and 
 \[
 2(1-\gamma_1)=(\H-\gamma_1S_2)
 \cdot \Gamma_1\geq (\gamma-\gamma_1)\Gamma\cdot\Gamma_1.
 \] 
It is easy to see that $\Gamma^Z$ and $\Gamma^Z_1$ intersect in at
least 2 nonsingular points on $S^Z_1$, and from this conclude that, in
this case also, $\gamma\leq 1$ (the details are very similar to
Cases~(a) and (b) and left to the reader).

It remains to show that Cases~(d), (e) can not occur, that is, $C =
\Pi\cap X$ is always reduced when $X$ is general.

\subparagraph*{\textsc{Claim}} If $X$ is general, every plane section 
of $X$ is reduced.

This is a fairly easy exercise. In coordinates $X$ is given by
 \[
 x_0^2 x_1x_2 + x_0 a_3 + b_4 =0
 \]
Where $a_3=a_3(x_1, ..., x_4)$ and $b_4=b_4(x_1, ..., x_4)$ are a
homogeneous cubic and quartic not involving $x_0$. The singular point 
$P \in X$ is of course the coordinate point $(1,0,0,0,0)$. Consider
the projection $\pi : X \dasharrow \P^3$ on $\P^3$
with homogeneous coordinates $x_1, ... , x_4$; it is a generically 2-to-1
map, which is to say that the equation of $X$ is
\emph{quadratic} in the variable $x_0$. Now $\Pi =\pi^{-1} \ell$ for a
unique line $\ell \subset \P^3$ and it is almost immediate that the 
hyperplane section $\Pi \cap X$ is nonreduced if, and only if,
either one of the following happens:
\begin{enumerate}
 \renewcommand{\labelenumi}{(\alph{enumi})}
 \item The line $\ell$ is contained in the discriminant surface
   $x_1x_2b_4 -a_3^2 =0$. It is very easy to see that, for a general
   choice of $a$, $b$, this surface contains no lines.
 \item The line $\ell$ is contained in the plane $x_1 =x_2$ \emph{and} 
   $a_3$, $b_4$, when restricted to $\ell$, both have a double root at 
   $x_1=x_2=0$. In any event, this is ruled out by condition~$2.2(b)$.
\end{enumerate}
\end{proof}

\section{Centers on $Y_{3,4}$} \label{sec:x34}

In this section we study maximal centers on $Y=Y_{3,4}$. We show
first that no curve on $Y$ can be a maximal center,
Theorem~\ref{thm:curve34}, then that a
nonsingular point can not be a maximal center,
Theorem~\ref{thm:point34}. 

\begin{thm} No curve on $Y$ can be a maximal center. \label{thm:curve34}
\end{thm}

\begin{proof} We can choose weighted projective coordinates 
$(x_1,x_2,x_3,x_4,y_1,y_2)$ such that the equations of $Y$ are as follows: 
 \[
Y:
\left\{\begin{array}{ll}
y_1y_2+b_4(x_1,x_2,x_3,x_4) & = 0 \\ 
y_1x_1+y_2x_2+a_3( x_1,x_2,x_3,x_4) & = 0
\end{array}
\right.
 \]
To understand the proof, it helps to know some explicit features of
the geometry of $Y$. To begin with, $Y$ is nonsingular apart from two 
$\mathbb{Z} /2\mathbb{Z}$-points $q_1$, $q_2$ at $(0,0,0,0,1,0)$ and
$(0,0,0,0,0,1)$. Denote $\rho_i: Y \dasharrow \P(1^4, 2)$ the
projection from $q_i \in Y$; it can be useful to know that the image 
of $\rho_1$, for example, is the hypersurface $x_2 y_2^2+a_3y_2
-x_1b_4=0$, as can be readily calculated eliminating the variable
$y_1$ from the equations of $Y$. Also, denote
$\pi:\P(1,1,1,1,2,2)\flip \P^3$ the projection on the coordinates of
degree $1$. 

The curves of $Y$, contracted by $\rho_1$, are the 12 curves $\ell$ with 
$\deg \O_\ell(1)=1/2$ given by $x_1=a_3=b_4=0$. Similarly, the
curves of $Y$, contracted by $\rho_2$, are the 12 curves $\ell$ with 
$\deg \O_\ell(1)=1/2$ given by $x_2=a_3=b_4=0$.
Finally, the curves contracted by $\pi$ are the 24 just mentioned,
plus the 3 curves $C$ with $\deg \O_C(1)=1$ given 
by $x_1=x_2=a_3=0$; under the generality assumption~$2.2(b)$ these are 
irreducible.

Assume that the curve $\Gamma = C_X(E)$ is the center of a maximal singularity 
$E$ of a mobile linear system $\H\subset |\O(n)|$.
By Corollary~\ref{cor:Kaw}, $\Gamma$ is contained in the nonsingular
locus of $X$. Denote $d=\deg \O_\Gamma(1)$, $m= \mult_{\Gamma}\H>n$. 

Choosing general members $H_1$, $H_2$ of $\H$ and 
intersecting with a general hyperplane section $S$ we obtain
 \[
 3n^2=H_1\cdot H_2\cdot S\geq m^2 d.
 \]
This implies that $d \leq 2$. We treat the two cases $d=2$, $d=1$ separately.

\paragraph*{\textsc{Case~$d=2$}} Here $\pi(\Gamma)$ is either a
line or a conic in $\P^3$. In either case $\Gamma$ is a
nonsingular rational curve and $\Gamma$ is defined scheme
theoretically by the quartics (with the natural
embedding of $\P(1^4,2^2)$ in $\P^{11}$ the curve $\Gamma$ is a normal 
quartic). Let $S\in |I_{\Gamma, Y}(4)|$ be a general member; write
as usual 
 \[
 A=\frac{1}{n} \H_{|S}=L+\gamma\Gamma,
 \]
(with $L=(1/n)\L$ nef...). We easily calculate on $S$ that $\Gamma^2=-8$, and
 \[
 0\leq L^2=12-4\gamma -8\gamma^2.
 \]
This implies that $\gamma\leq 1$ and finishes this case.

\paragraph*{\textsc{Case~$d=1$}} Here $\pi(\Gamma)\subset \P^3$ is a line,
$\Gamma$ is a nonsingular rational curve. Denote $S_1$, $S_2$ two
general members of the pencil $|I_{\Gamma, Y}(1)|$, $C =\bs |S_1 \cap
S_2|$ the base locus. Denoting $\Pi=\pi^{-1} \pi \Gamma \cong \P(1,1,2,2)$, 
we can also say that $C= \Pi\cap Y$.

In the end we will show that $C$ is reduced; for now let us assume it.
We restrict to $S_1$ and write
\begin{eqnarray*}
A =(1/n)\H_{|S_1}=& L+\gamma\Gamma+\sum_{i=1}^r \gamma_i\Gamma_i \\
       S_{2|S_1} =& C=\Gamma+\sum \Gamma_i
\end{eqnarray*}
When $C$ is reduced, it is clear that if $W$ is a component of $C$, 
$\mult_W \H=\mult_W\H_{|S_1}$; in particular $\mult_\Gamma \H
=\mult_\Gamma \H_{|S_1}=n\gamma$, and also $\mult_{\Gamma_i} \H =
n\gamma_i$. Because $\Gamma$ is a \emph{maximal} singularity, $\gamma 
\geq \gamma_i$ for all $i$, hence possibly after relabelling
components of $C$, we can assume that $\gamma\geq \gamma_r\geq ...\geq
\gamma_1$. 
Consider the effective $\Q$-divisor
 \[
 M=(A-\gamma_1S_2)_{|S_1}=L+(\gamma-\gamma_1)\Gamma
 +\sum_{i>1}(\gamma_i-\gamma_1)\Gamma_i.
 \]
We calculate the intersection product, on $S_1$, with $\Gamma_1$:
 \[
M\cdot \Gamma_1 = (1-\gamma_1)\deg \O_{\Gamma_1}(1) \geq (\gamma
-\gamma_1) \Gamma \cdot \Gamma_1
 \]
Note that here $1/2 \leq \deg \O_{\Gamma_1}(1) \leq 2$ is a half-integer. 
It is completely elementary to check that $\Gamma\cdot\Gamma_1\geq
\deg \O_{\Gamma_1}(1)$ (in doing this, it helps to note that $S_1$ is
nonsingular outside $q_1$, $q_2$). Together with the last displayed
equation this implies that $\gamma \leq 1$ and finishes the proof. 
It remains to show that $C$ is reduced.

\paragraph*{\textsc{Claim}} If $X$ is general $C$ is reduced.

This is a fairly easy exercise; the situation corresponds exactly to
the quartic $X_4$ as treated in the proof of Theorem~\ref{nocurves}.
In short, it is easy to see that $C= \Pi \cap Y$ is nonreduced
if, and only if, either one of the following happens:
\begin{enumerate}
 \renewcommand{\labelenumi}{(\alph{enumi})}
 \item The line $\pi \Gamma$ is contained in the discriminant surface
   $x_1x_2b_4 -a_3^2 =0$. It is clear that, for a general
   choice of $a$, $b$, this surface contains no lines.
 \item The line $\pi \Gamma$ is contained in the plane $x_1 =x_2$ \emph{and} 
   $a_3$, $b_4$, when restricted to $\pi \Gamma$, both have a double root at 
   $x_1=x_2=0$. In any even this possibility is certainly ruled out by
   condition $2.2(b)$. 
\end{enumerate}

\end{proof}

\begin{thm} \label{thm:point34} Let $x \in Y$ be a nonsingular point. 
Then $x$ is not a maximal center.
\end{thm}

\begin{proof} Let $x\in Y$ be a nonsingular point, $B= \bs
  |I_x(1)|$. If $\dim B=0$ consider a general element
  $S\in|I_x(1)|$. 
Then $\H_{|S}$ is mobile and
 \[
 3=\frac{1}{n^2}\H^2\cdot S<4.
 \]
This is enough, by Theorem~\ref{thm_2ds}, to conclude that $x$ is not a  
center. 

Let us now worry about the case $\dim B=1$; this can only happen when 
$x \in B$ is a curve contracted by $\pi$, and as we have already noted
at the start of the proof of Theorem~\ref{thm:curve34}, it is a
consequence of the generality assumption~$2.2(b)$ that $B$ is irreducible. 
If $\deg \O_B(1)=1$ then write $\H_{|S}=\L+mB$ where $\L$ is the mobile part;
dividing by $n$ be obtain
 \[
 A=\frac{1}{n}\H_{|S}=L+cB
 \]
where $L=(1/n)\L$, and $c=m/n$. Note that $B$
is a rational curve on $S$ passing trough 2 simple double points.
Therefore
$(B\cdot B)_S=-1$ and,
computing the self intersection of $L$, we get
 \[
 L^2=3-c^2-2c\leq 4(1-c).
 \]
Again by Theorem~\ref{thm_2ds} we exclude $x$ as a center.

If $\deg \O_B(1)B=1/2$ then the above arguments on the surface $S$
are not enough to exclude it as a maximal center. This time we need to 
consider the linear system $\Delta=|I_x^{\otimes 2}(2)|$.

It is easy to check that:

\textsc{Claim} If $D\in\Delta$ is general, then
\begin{enumerate}
 \renewcommand{\labelenumi}{(\alph{enumi})}
 \item $D$ has a simple double point at $x$,
 \item $D$ is nonsingular along $B\setminus(\Sing(Y)\cup\{x\})$,
 \item $D$ has a singularity of type $1/4(1,-1)$ at $B\cap \Sing(Y)$.
\end{enumerate}

Let $\nu:Y^{\prime}\to Y$ the blow up of $x$ with exceptional divisor
$E$. Write $\nu^*D=D^{\prime}+bE$ and $F=E_{|D^{\prime}}$.
$D$ has a double point at $x$ thus $F\subset E$ is a conic.
By Shokurov connectedness there is a line $\ell\subset E$
such that 
 \[
 \ell \subset LC\Bigl(Y^{\prime},
 K_{Y^\prime}+\frac{1}{n}\H^{\prime}+\Bigl(\frac{b}{n}-1\Bigr)E
 +D^{\prime}\Bigr).
 \]
Therefore, for the generic $D$, by inversion of adjunction
$K_D^{\prime}+((1/n)\H^{\prime}+(b/n+1)E)_{|D^{\prime}}$ is not
LC at two distinct points $p_1$ and $p_2$.
We want to use this fact to derive a numerical constraint
on $\H$. To do this let us first compute the 
intersection matrix on $D^{\prime}$:
 \[
 (F\cdot F)_{D^{\prime}}=-2, \hspace{.2cm} 
 (B^{\prime}\cdot B^{\prime})_{D^{\prime}}=-\frac74,
 \hspace{.2cm} (F\cdot B^{\prime})_{D^{\prime}}=1,
 \]
(The only nontrivial product is $(B^{\prime}\cdot B^{\prime})_{D^{\prime}}
=-2-1/2+3/4$, by the adjunction formula with correction coming from
the different). Write
 \[
 A=\Bigl(\nu^* \frac{1}{n}\H \Bigr)_{|D^{\prime}}=
 \frac{1}{n} \L+\beta F+\alpha B^{\prime},
 \]
where $\L$ is the mobile part. We have
 \begin{equation}\label{vero}
 (\L/n)^2=6-2\beta^2-\frac74\alpha^2-\alpha+2\beta\alpha.
 \end{equation}
To find a lower bound for $(\L/n)^2$ recall that
$K_D^{\prime}+(\H^{\prime}/n+(b/n-1)E)_{|D^{\prime}}$ is not log canonical at 
$p_1$ and $p_2$ therefore by Theorem~\ref{thm_2ds}, we always have 
 \[
 (\L/n)^2>4(2-\beta)+4(2-\beta)(1-\alpha)=16-8\beta-8\alpha+4\beta\alpha.
 \]
Combining with equation~\ref{vero} yields 
 \[
 0>2\beta^2-2\beta(4-\alpha)-7\alpha+\frac74\alpha^2+10,
 \]
and the discriminant of this quadratic equation with respect to $\beta$ is
 \begin{multline*}
 \Delta/4=16+\alpha^2-8\alpha+14\alpha-\frac72\alpha^2-20 \\
 =-\frac52\alpha^2+6\alpha-4
 =-\frac52\bigl(\alpha-\frac65\bigr)^2-\frac45<0.
 \end{multline*}
This inequality shows that $x$ cannot be a center of maximal singularities
and concludes the proof of the Theorem.

\end{proof}

\end{document}